\documentclass[12pt]{article}
\usepackage{color}
\usepackage{epsfig,subfigure,epstopdf}
\usepackage{graphicx}
\usepackage{natbib}
\usepackage{amsmath,bm}
\usepackage{amssymb}
\usepackage{amsthm}
\graphicspath{{figures/}}

\def\red{\color{red}}

\begin{document}
\title{On the Weak Convergence and Central Limit Theorem of Blurring
and Nonblurring Processes with Application to Robust Location
Estimation}

\author{Ting-Li Chen$^a$\footnote{Corresponding author:
tlchen@stat.sinica.edu.tw},~
Hironori Fujisawa$^b$,~ Su-Yun Huang$^a$,~ Chii-Ruey Hwang$^c$\\[3ex]
$^a$ Institute of Statistical Science, Academia Sinica, Taiwan\\
$^b$ Institute of Statistical Mathematics, Japan\\
$^c$ Institute of Mathematics, Academia Sinica, Taiwan }

\date{\today} \maketitle

\newtheoremstyle{dotless}{}{}{\itshape}{}{\bfseries}{}{ }{}{}

\renewcommand{\labelenumi}{(\roman{enumi})}
\newtheorem{example}{Example}
\newtheorem{definition}{Definition}
\newtheorem{theorem}{Theorem}
\newtheorem{lemma}{Lemma}
\newtheorem{corollary}{Corollary}
\newtheorem{remark}{Remark}

\def\1{{(1)}} 
\def\t{{(t)}}\def\tplus{{(t+1)}}\def\tn{{[t]}}\def\tnplus{{[t+1]}}
\def\it{{-(t)}}\def\-itn{{-[t]}}\def\itplus{{-(t+1)}}\def\itnplus{{-[t+1]}}
\def\s{{(s)}}\def\splus{{(s+1)}}\def\sn{{[s]}}\def\snplus{{[s+1]}}
\def\is{{-(s)}}\def\-isn{{-[s]}}\def\isplus{{-(s+1)}}\def\isnplus{{-[s+1]}}

\def\T{{(\tau)}} \def\Tn{{[\tau]}}
\def\Real{{\mathbb R}}\def\real{{\mathbb R}}
\def\th{{\rm th}}
\def\cov{{\rm Cov}}
\def\hat{\widehat}
\def\tilde{\widetilde}

\def\red{\color{red}}
\def\blue{\color{blue}}
\def\hiro{\color[rgb]{0,0.5,0}}

\begin{abstract}
This article studies the weak convergence and associated Central
Limit Theorem for blurring and nonblurring processes. Then, they are
applied to the estimation of location parameter. Simulation studies
show that the location estimation based on the convergence point of
blurring process is more robust and often more efficient than that of
nonblurring process.

\noindent {\bf Keywords}: Weak convergence, Central limit theorem, Blurring process, Robust estimation.
\end{abstract}

\section{Introduction}

In this article we consider two types of processes arisen from
mean-shift algorithms \citep{Fukunaga,ChengY}. Starting with $n$
points $\{x_{i,n}\}_{i=1}^n$ as initials, the nonblurring type
process is given by
\begin{eqnarray}
\label{eq:updatnb}
x_{i,n}^\tnplus
 =\frac{\int x\,w(x-x_{i,n}^{[t]})dF_n(x)} {\int w(x- x_{i,n}^{[t]})dF_n(x)}
 = \sum_{j=1}^N  \displaystyle
\frac{w(x_{j,n} - x_{i,n}^\tn )}{\sum_{\ell=1}^N
w(x_{\ell,n} -x_{i,n}^\tn)}\, x_{j,n},
\end{eqnarray}
where $x^{[0]}_{i,n}=x_{i,n}$ and $F_n$ is the empirical
distribution function based on the initial points
$\{x_{i,n}\}_{i=1}^n$. This process~(\ref{eq:updatnb}) consists of
$n$ simultaneous updating paths, wherein each path starts from one
initial. Another type of updating process, called the blurring type,
is considered by replacing $F_n$ with the iteratively updated
empirical distribution $F_n^\t$ based on updated points
$\{x_{i,n}^\t\}_{i=1}^n$, in addition to the above idea of weighted
scores for updating:
\begin{eqnarray}
\label{eq:updatbl}
x_{i,n}^{(t+1)}
 =\frac{\int x\,w(x- x_{i,n}^{(t)})dF_n^{(t)}(x)} {\int w(x- x_{i,n}^{(t)})dF_n^{(t)}(x)}
 = \sum_{j=1}^N  \displaystyle
\frac{w(x_{j,n}^{(t)}-x_{i,n}^{(t)}) } {\sum_{\ell=1}^N
w(x_{\ell,n}^{(t)}-x_{i,n}^{(t)})}\, x_{j,n}^{(t)},
\end{eqnarray}
where $x^{(0)}_{i,n}=x_{i,n}$. Same as the nonblurring process, the
blurring process~(\ref{eq:updatbl}) starts with $n$ initials
$\{x_{i,n}\}_{i=1}^n$, and then it goes through a simultaneous
updating at each iteration. The key difference from the nonblurring
process is that, this process~(\ref{eq:updatbl}) takes weighted
average according to the updated empirical distribution $F_n^\t$,
while the nonblurring process takes weighted average with respect to
the initial empirical distribution $F_n$. That is, at each iteration
in the blurring process, not just the weighted centers are updated
from $\{x_{i,n}^\t\}_{i=1}^n$ to $\{x_{i,n}^\tplus\}_{i=1}^n$, the
empirical distribution is also updated from $F_n^\t$ to
$F_n^\tplus$.

The blurring process was developed and named SUP (self-updating
process in Chen and Shiu, 2007) and was recently applied to cryo-em
image clustering (Chen et al., 2014). It is also known as the
blurring type mean-shift algorithm (Cheng, 1995; Comaniciu and Meer,
2002). Blurring mean-shift can be viewed as a homogeneous
self-updating process. Algorithm convergence and location estimation
consistency of the blurring and nonblurring processes were discussed
in \citet{ChengY}, \citet{Comaniciu}, \citet{Li}, \citet{Chen2}, and
\citet{Ghassabeh}. In this article, we study their weak convergence
and associated Central Limit Theorem. Due to the complicated
dependent structure of random variables in $\{x_{i,n}^\tn\}_{i=1}^n$
and $\{x_{i,n}^\t\}_{i=1}^n$, the study of their asymptotic behavior
becomes challenging.

The convergence point of the blurring or the nonblurring process can
be used for location estimation, which is one of the most basic and
commonly used tasks in statistical analysis as well as in computer
vision. It is well-known that the sample mean is not a robust
location estimator and it is sensitive to outliers and data
contamination. To reduce the influence from deviant data, there is a
wide class of robust M-estimators in statistics literature using
weighted scores (Hampel et al., 1986; Huber, 2009; van de Geer,
2000). Consider a weighted score equation for the mean $\mu$:
\begin{eqnarray}\label{est_eq}
\sum_{i=1}^n w(x_{i,n}-\mu) (x_{i,n}-\mu) = 0,
\end{eqnarray}
where $w(x)$ is a symmetric weight function. The weighted mean that
satisfies the estimating equation~(\ref{est_eq}) can be shown to
take the following form
\begin{eqnarray}\label{mu_eq}
\mu = \frac{\sum_{i=1}^n x_{i,n} w(x_{i,n}-\mu) }{\sum_{i=1}^n w(x_{i,n}-\mu)}
 = \frac{\int x w(x-\mu) dF_n(x)}{\int w(x-\mu) dF_n(x)}.
\end{eqnarray}
This estimator~(\ref{mu_eq}) can be obtained by the fixed-point
iteration algorithm at convergence, where the iterative update is
given by
\begin{eqnarray}
\mu^\tnplus = \frac{\int x w(x-\mu^\tn) dF_n(x)}{\int w(x-\mu^\tn) dF_n(x)}, \qquad t=0,1,2,\ldots.
\label{eq:score2}
\end{eqnarray}
with an appropriate starting initial $\mu^{[0]}$. Here, we consider
a simple change of the updating process by starting with $n$ data
points $\{x_{i,n}\}_{i=1}^n$ as initials and by replacing $\mu^\tn$
with $x_{i,n}^\tn$ in~(\ref{eq:score2}). It then leads to the
nonblurring process given in~(\ref{eq:updatnb}):
\[
x_{i,n}^\tnplus
 =\frac{\int x\,w(x-x_{i,n}^{[t]})dF_n(x)} {\int w(x- x_{i,n}^{[t]})dF_n(x)}
 = \sum_{j=1}^N  \displaystyle
\frac{w(x_{j,n} - x_{i,n}^\tn )}{\sum_{\ell=1}^N
w(x_{\ell,n} -x_{i,n}^\tn)}\, x_{j,n}.
\]
By replacing $F_n$ with $F_n^\t$, we have the blurring process given
in~(\ref{eq:updatbl}):
\[
x_{i,n}^{(t+1)}
 =\frac{\int x\,w(x- x_{i,n}^{(t)})dF_n^{(t)}(x)} {\int w(x- x_{i,n}^{(t)})dF_n^{(t)}(x)}
 = \sum_{j=1}^N  \displaystyle
\frac{w(x_{j,n}^{(t)}-x_{i,n}^{(t)}) } {\sum_{\ell=1}^N
w(x_{\ell,n}^{(t)}-x_{i,n}^{(t)})}\, x_{j,n}^{(t)}.
\]
The iterative updating process based on either (\ref{eq:updatnb}) or
(\ref{eq:score2}) has been adopted for robust mean estimation (Field
and Smith, 1994; Fujisawa and Eguchi, 2008; Maronna, 1976; Windham,
1995; among others) and robust clustering (Notsu et al., 2014). It
is also known as the nonblurring type mean-shift algorithm. On the
other hand, robust estimation based on blurring approach is rather
rare in the literature. Here we strongly recommend it as an
alternative choice. From our simulation studies in
Section~\ref{sec:simulation}, the blurring type algorithm is often
more robust with smaller mean square error. Thus, the blurring type
algorithm deserves more attention and further exploration.

The contribution of this article is twofold. First, we derive
theoretical properties of the blurring and nonblurring processes
including their weak convergence to a Brownian bridge-like process
and associated Central Limit Theorem. These theoretical results are
presented in Section~\ref{sec:main}, with all technical proofs being
placed in the Appendix. Second, we apply the derived Central Limit
Theorem to location estimation. Simulation studies comparing
location estimation based on using blurring and nonblurring
processes are presented in Section~\ref{sec:simulation}. Our
simulation results suggest that the blurring type algorithm is often
more robust than the existent nonblurring type algorithm for robust
M-estimation.

\section{Main Results}\label{sec:main}

Let $\{x_{i,n}\in\real: i=1,\dots,n\}$, $n\in{\mathbb N}$, be a
triangular array of random variables. Assume the following
conditions.
\begin{itemize}\itemsep=0pt
\item[C1.]
The underlying distribution $F$ has a continuous probability
density function $f(x)$, which is symmetric about its mean
$\mu$.
\item[C2.]
The weight function $w(x)$ is a probability density function. It
is log-concave and symmetric about 0. (This condition implies
that $w(x)>0$ for all $x\in\real$ and that $w(x)$ is unimodal
and non-constant.)
\item[C3.]
For simplicity but without loss of generality, assume $\mu=0$.
\end{itemize}

\subsection{Weak convergence of blurring process}

Let $x^{(0)}:=x$ and
\begin{eqnarray*}
x^\tplus &:=&\eta^{(t+1)}(x^{(t)}) =\eta^{(t+1)}\circ\cdots\circ \eta^{(1)}(x),
\quad t=0,1,2,\ldots,
\end{eqnarray*}
where $\eta^\tplus$ is a blurring transformation given by
\begin{equation}\label{eta}
\eta^{(t+1)}(x):=\frac {\int y \cdot w(y-x) dF^{(t)}(y)} {\int w(y-x) dF^{(t)}(y)},
\end{equation}
where $\{F^\t\}_{t=0,1,\ldots}$ is the cumulative distribution
functions of $X^\t$ with $X^{(0)} \sim F$.  Note that the blurring
transformation $\eta^\tplus$ shifts $x$ toward a mode by an amount
depending on $F^\t$,
\[\frac {\int (y-x) \cdot w(y-x)dF^{(t)}(y)} {\int w(y-x) dF^{(t)}(y)}.\]
Let  $\eta_n^\tplus$ be the empirical blurring transformation based
on $F_n^\t$, which is the empirical cumulative distribution of
$\{x_{i,n}^\t\}_{i=1}^n$. In Theorem~\ref{thm:main} below, it is
shown $F_n^\t(x) \to F^\t(x)$ almost surely for each $x$, as
$n\to\infty$. The blurring process, in empirical level and in
population level, can be expressed as
\begin{eqnarray*}
x_{i,n}^{(t+1)} &=& \eta_n^{(t+1)}(x_{i,n}^{(t)}) =\eta_n^{(t+1)}\circ\cdots\circ \eta_n^{(1)}(x_i^{(0)})
\quad\mbox{(empirical blurring)},\\
x^{(t+1)} &=&\eta^{(t+1)}(x^{(t)})=\eta^{(t+1)}\circ\cdots\circ
\eta^{(1)}(x^{(0)}) \quad\mbox{(population blurring)}.
\end{eqnarray*}

At $t=0$, it is known that $F_n(x) \to F(x)$ almost surely for each
$x$, as $n \to \infty$. Furthermore, by Donsker's Theorem, the
sequence
$$Z_n(x):=\sqrt{n}\left(F_n(x)-F(x)\right)$$
converges in distribution to a Gaussian process with zero mean and
covariance given by
\[\cov\left(Z_n(x),Z_n(y)\right) =\min\{F(x),F(y) \} -F(x)F(y).\]
Let this convergence in distribution be denoted by
\[
Z_n(x):=\sqrt{n}\left(F_n(x)-F(x)\right) \rightsquigarrow
B(F(x)),
\]
where $B$ is the standard Brownian bridge. In this article we
establish the weak convergence for the empirical process of
cumulative distribution function for each iteration. We will show
the weak convergence by mathematical induction. The weak convergence
is true at $t=0$ by Donsker's Theorem. Next, by assuming that
$F_n^\s(x) \to F^\s(x)$ almost surely and that $
\sqrt{n}\left(F_n^\s(x)-F^\s(x)\right) \rightsquigarrow B^\s(F(x))$
for some Brownian bridge like process $B^\s$, we show that claimed
statements hold for $t=s+1$. Because of the complicated dependent
structure in $\{x_{i,n}^\t\}$, the almost sure convergence for
$F_n^\t(x)$ and the weak convergence for $Z_n^\t(x)$ become
difficult, where
$$Z_n^\t(x):=\sqrt{n}\left( F^\t_n(x) - F^\t(x) \right).$$
We first establish the connection between the empirical process of
cumulative distribution functions of two consecutive iterations.
Then we prove that this connection is a continuous mapping under the
Skorokhod topology.

Before establishing the main theorem we derive a few technical
lemmas first. Lemma~\ref{lemma:order}, with the proof given in
Appendix~\ref{appendix:lemma1}, shows that $\eta^{ \t}$ is a
one-to-one transformation, which implies that the data orders do not
change during the blurring process at each update. This phenomenon
is important when we calculate the empirical cumulative distribution
function of the current iteration based on the process of previous
iteration.

\begin{lemma}\label{lemma:order}
Assume conditions C1-C3. We have $0<\eta^\tplus(a)<\eta^\tplus(b)$
for $0<a<b$.
\end{lemma}

Next in Lemma~\ref{lemma2}, we derive the asymptotic behavior of
$\eta_n^\tplus$, which immediately implies that $F_n^\tplus(x) \to
F^\tplus(x)$ almost surely.
\begin{lemma}\label{lemma2} For any $x$,
$\lim_{n \to \infty} \eta_n^\tplus (x) = \eta^\tplus(x)$ a.s. for
$t=0,1,\dots$
\end{lemma}
The proof is given in Appendix~\ref{appendix:lemma2}.

To show the weak convergence of $Z_n^\tplus(x)$, we need a tighter
estimate of $\eta_n^\tplus-\eta^\tplus$, which is presented in the
following lemma with the proof given in
Appendix~\ref{appendix:lemma3}. While Lemma~\ref{lemma2} shows
$\eta_n^\tplus(x)-\eta^\tplus(x)=o_p(1)$, Lemma~\ref{lemma3} implies
that $\eta_n^\tplus -\eta^\tplus(x)=O(1/\sqrt{n})$.
\begin{lemma}\label{lemma3}
For $t=0,1,2,\dots$,
\begin{equation}
\sqrt{n}\left(\eta_n^\tplus(x)-\eta^\tplus(x)\right)
=\frac {\int \left(y-\eta_n^\tplus(x)\right) w(y-x)\,
 d Z_n^\t(y)}{\int  w(y-x) dF^\t(y)}.\label{eq:wc_g}
\end{equation}
\end{lemma}

Let $\xi^\t$ and $\xi_n^\t$ be the inverse function of $\eta^\t$ and
$\eta_n^\t$, respectively. Then, we have
\begin{eqnarray*}
F^\tplus(x) = P\left\{ X^\tplus \le x\right\} = P\left\{ \eta^\tplus(X^\t)\le x\right\} =
F^\t({\xi^\tplus}(x)).
\end{eqnarray*}
Using this formula and the Taylor expansion, we establish the
connection between $Z_n^\tplus(x)$ and $Z_n^\t(x)$. The result is
presented in the following lemma with the proof given in
Appendix~\ref{appendix:lemma4}.

\begin{lemma}\label{lemma4}
The process $Z_n^\tplus(x)$ can be expressed as a function of the
previous process $Z_n^\t(x)$ up to an additive $o_p(1)$-term.
Precisely,
\begin{eqnarray}\label{eq:wc_it}
Z_n^\tplus(x) = \int_y \left\{K^\tplus(x,y)+o_p(1)\right\}
  dZ_n^\t(y),
\end{eqnarray}
where
\begin{equation}\label{eq:Gt}
K^\tplus(x,y)= -\frac{f^\tplus(x)\cdot(y-x) w(y-\xi^\tplus(x))}{ {\int  w\left(y-\xi^\tplus(x)\right) dF^\t(y)}}
+1_{\{y\leq \xi^\tplus(x)\}}.
\end{equation}
\end{lemma}
With the assumption of the weak convergence of $Z_n^\t(x)$ to a
Brownian bridge-like process, the variance of $\int_y o_p(1)
dZ_n^\t(y)$ will also be $o_p(1)$.

Lemma~\ref{lemma:sk_conv} below states that the mapping
(\ref{eq:wc_it}), ignoring $o_p(1)$, is continuous. The continuity
is stated in Lemma~\ref{lemma:sk_conv} with the proof given in
Appendix~\ref{appendix:lemma5}.

\begin{lemma}\label{lemma:sk_conv}
Let ${\cal D}={\cal D}[0,1]$ be the space of real-valued functions
on $[0,1]$ that are right-continuous and have left-hand limits.
Define ${\cal L}^\tplus$ as a mapping from ${\cal D}$ to ${\cal D}$,
such that
\begin{equation}\label{cal_L}
({\cal L}^\tplus  W) (u)=
\int_{v=0}^1 K^\tplus \left(F^\it(u),F^\it(v)\right) \, dW(v) \quad \forall W \in {\cal D},
\end{equation}
where $F^\it(u)$ denotes the inverse function of $F^\t(x)$. Then,
${\cal L}^\tplus$ is continuous under the Skorokhod topology.
\end{lemma}

Let $B_n^\t(u):= Z_n^\t(F^\it(u))$. Then $B_n^\t(u) \in {\cal D}$,
and
\[
B_n^\tplus(u) = \int_{v=0}^1 \left\{K^\tplus(F^\it(u),F^\it(v))+o_p(1)\right\}
  dB_n^\t(v).
\]
Lemma \ref{lemma:sk_conv} shows that the mapping is asymptotically
continuous. Therefore, $B_n^\tplus$ has the same weak convergence
property as $B_n^\t$. The result is summarized in
Theorem~\ref{thm:main}.

\begin{theorem}\label{thm:main}
Assume conditions C1-C3. We have
\begin{enumerate}\itemsep=0pt
\item
For each $x$, $F_n^\t(x) \to F^\t(x)$ almost surely, as
$n\to\infty$;
\item
$\sqrt{n}\left(F_n^{(t)}(x)-F^{(t)}(x)\right) \rightsquigarrow
B^{(t)}(F^\t(x))$, where $B^{(0)}$ is the standard Brownian
bridge on $[0,1]$, $B^\t(u)=\int_{v=0}^1 H^\t(u,v) dB^{(0)}(v)$,
$H^\1 (u,z):= K^\1(F^{-1}(u),F^{-1}(z))$, and
\[
H^\tplus(u,z):= \int_{v=0}^1 K^\tplus(F^\it(u),F^\it(v))
\frac {\partial H^\t(v,z) }{\partial v}\, dv.
\]
\end{enumerate}
\end{theorem}
Proof for Theorem~\ref{thm:main} is in Appendix~\ref{appendix:thm1}.


\begin{corollary}
$B^\t(0)=B^\t(1)=0$ and ${\cov\left\{ B^\t(u), B^\t(v)\right\}=
\int_0^1 H^\t(u,z) H^\t(v,z)dz}$.
\end{corollary}

\subsection{Weak convergence of nonblurring process}
In the nonblurring process we use almost the same notation as in the
blurring process except for the superscript $[t]$. A similar result
for the nonblurring process is stated in the following theorem with
its proof given in Appendix \ref{appendix:thm2}.

\begin{theorem}\label{thm:main2}
Assume conditions C1-C3. We have
\begin{enumerate}
\item
For each $x$, $F_n^\tn(x) \to F^\tn(x)$ almost surely, as
$n\to\infty$;
\item
$\sqrt{n}\left(F_n^\tn(x)-F^\tn(x)\right) \rightsquigarrow
B^\tn(F(x))$, where $B^{(0)}$ is the standard Brownian bridge on
$[0,1]$, $B^\tn(u)=\int_{0}^1 H^\tn(u,v) dB^{(0)}(v)$ with
\begin{eqnarray*}
H^{[1]}(u,v)~~\, &:=&K^{[1]}(F^{-1}(u),F^{-1}(v))+
  1_{\left\{\eta^{[1]}(F^{-1}(v)) \leq F^{-1}(u)\right\}}\, ,\\
H^\tnplus(u,v) &:= & K^\tnplus(F^{-1}(u)),F^{-1}(v))
 +   H^\tn \left(F(\xi^\tnplus(F^{-1}(u))),F(\xi^\tnplus(F^{-1}(v)))\right).
\end{eqnarray*}
\end{enumerate}
\end{theorem}

\subsection{Central Limit Theorem}

Apply the weak convergence, we can have the Central Limit Theorem
for the sample mean of updated data. The result on the blurring case
is presented below with its proof given in Appendix
\ref{appendix:thm3}
\begin{theorem}\label{thm:clt_bl}
Let $S_n^\tplus :=\frac1{\sqrt n}\sum_{i=1}^n x_{i,n}^\tplus$. We
have
\[
S_n^\tplus \rightsquigarrow  \int\int \left(\frac{w(x-y)x}{\rho^\t(y)} +
\frac{w(y-x)y}{\rho^\t(x)}\right) dF^\t(x) dB^\t(F^\t(y)),
\]
where $\rho^\t(y)= E_{X \sim F^\t} \left(w(X-y)\right)$.
\end{theorem}

From Theorem \ref{thm:main},
\[
B^\t(u)=\int_{v=0}^1 H^\t(u,v) dB^{(0)}(v).
\]
Take the derivative, we have
\[
dB^\t(u)=\int_{v=0}^1 \frac {\partial H^\t(u,v)}{\partial u} dB^{(0)}(v) du.
\]
Now substitute $u=F^\t(y)$ in the above equation. We have
\[
dB^\t(F^\t(y))=\int_{v=0}^1 \frac {\partial H^\t(u,v)}{\partial u}\Big|_{u=F^\t(y)} dB^{(0)}(v) dF^\t(y).
\]
Therefore
\begin{eqnarray*}
&&\int\int \left(\frac{w(x-y)x}{\rho^\t(y)} +
\frac{w(y-x)y}{\rho^\t(x)}\right) dF^\t(x) dB^\t(F^\t(y))\\
&=&\int\int \left(\frac{w(x-y)x}{\rho^\t(y)} +
\frac{w(y-x)y}{\rho^\t(x)}\right) dF^\t(x) \int_{v=0}^1
\frac {\partial H^\t(u,v)}{\partial u}\Big|_{u=F^\t(y)}  dB^{(0)}(v) dF^\t(y).\\
&=& \int_{v=0}^1 \left[\int\int \left(\frac{w(x-y)x}{\rho^\t(y)} +
\frac{w(y-x)y}{\rho^\t(x)}\right)
\frac {\partial H^\t(u,v)}{\partial u}\Big|_{u=F^\t(y)}  dF^\t(x)dF^\t(y) \right]dB^{(0)}(v)
\end{eqnarray*}
This distribution has mean 0, and variance
\[
\int_{v=0}^1 \left[\int\int \left(\frac{w(x-y)x}{\rho^\t(y)} +
\frac{w(y-x)y}{\rho^\t(x)}\right)
\frac {\partial H^\t(u,v)}{\partial u}\Big|_{u=F^\t(y)}  dF^\t(x)dF^\t(y) \right]^2dv.
\]

The Central Limiting Theorem for the nonblurring case is similar,
and the proof is almost identical and thus is omitted. The update of
the nonblurring process is a weighted average over the original data
and the weights depend on the updated data at the previous
iteration.

\begin{theorem} Let $S_n^\tnplus :=\frac1{\sqrt n}\sum_{i=1}^n x_{i,n}^\tnplus$.
 We have
\begin{eqnarray*}
S_n^\tnplus &\rightsquigarrow & \int\int \frac{w(x-y)x}{\rho(y)}dF(x) dB^\tn(F^\tn(y))
+\int\int
\frac{w(y-x)y}{\rho(x)} dF^\tn(x) dB^{(0)}(F(y)),
\end{eqnarray*}
where $\rho(y)= E_{X\sim F} \left(w(X-y)\right)$.
\end{theorem}

\section{Application to robust location estimation with simulation studies}\label{sec:simulation}
The convergence point with either the blurring or the nonblurring
process is a reasonable robust estimation of the location parameter.
The consistency of the nonblurring process is proved
in~\citet{ChengY}, and that of the blurring process is proved
in~\citet{Chen2}. With the Central Limit Theorem provided in the
previous section, it is of our interest to compare the efficiency of
both processes. Theoretical comparison of asymptotic variance is
quite difficult even for the simplest case that both the sampling
distribution and the weight function are normal. In below we will
first explain the rationale behind the phenomenon that blurring is
more robust and often more efficient than nonblurring. Then, we will
show by simulation studies the behavior of asymptotic normality and
the mean square error comparison for both processes.

Without the update of empirical distribution i.e., by keeping the
initial empirical distribution $F_n$ throughout all iterations in
the nonblurring process, the effective weights in~(\ref{eq:updatnb})
are approaching
\[\frac{w(x-\mu^{[\tau]})}{\int w(x-\mu^{[\tau]}) dF_n(x)},\]
where $\tau$ is the last iteration step at convergence. With the
update of empirical distribution in the blurring process, the
effective weights in~(\ref{eq:updatbl}) are getting more and more
close to uniform. Since the weighted average is taken with respect
to newly updated centers, each of which is a weighted average of
previous updated centers, the contribution of each original data
point to the final estimation is relatively uniform for the blurring
process, while the contribution of each data point in the
nonblurring process is governed by $w(x)$. It is known that the
sample mean, which corresponds to a uniform weight, is the uniformly
minimum variance unbiased estimator for many distributions including
normal distribution. While both blurring and nonblurring estimators
are robust by reducing the contribution of outliers or data points
in heavy tails, the blurring estimator, which takes a relatively
uniformly weight, is expected to be more efficient than the
nonblurring estimator on processing the relatively reliable part of
information. Our simulation results presented in
Section~\ref{compare_var} also support this thinking.

\subsection{Asymptotic normality}

\begin{figure}[htb]
\includegraphics[width=6.5in,height=2.5in]{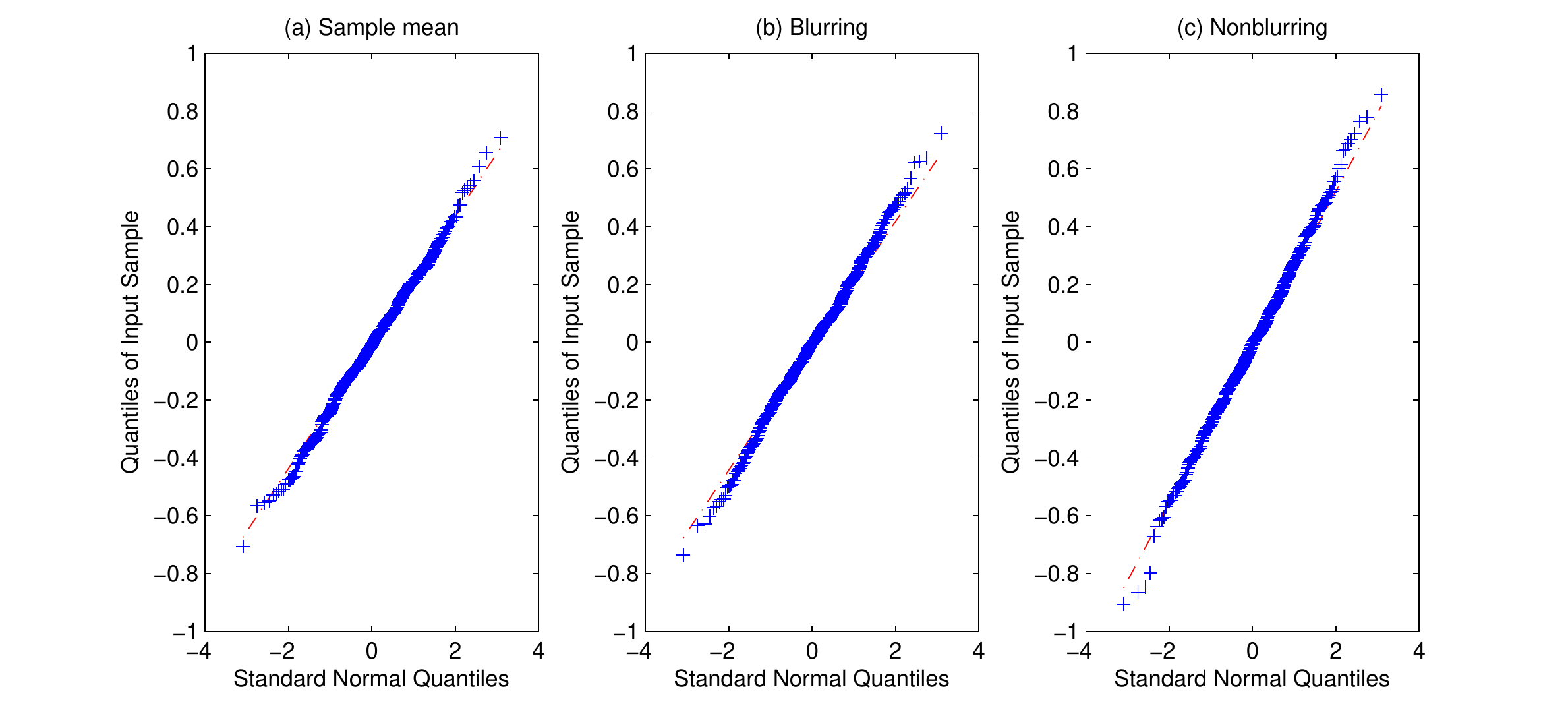}\\
\includegraphics[width=6.5in,height=2.5in]{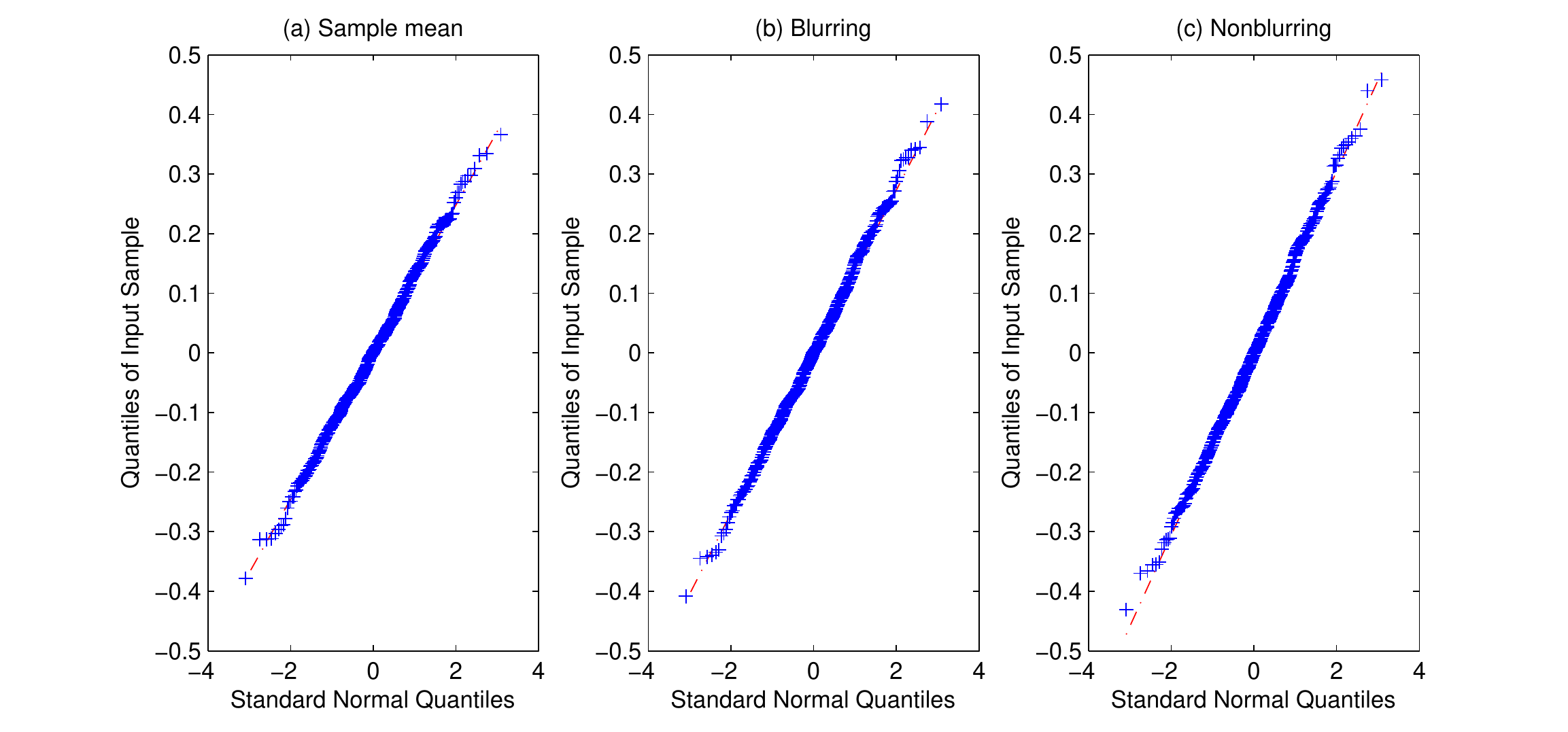}\\
\includegraphics[width=6.5in,height=2.5in]{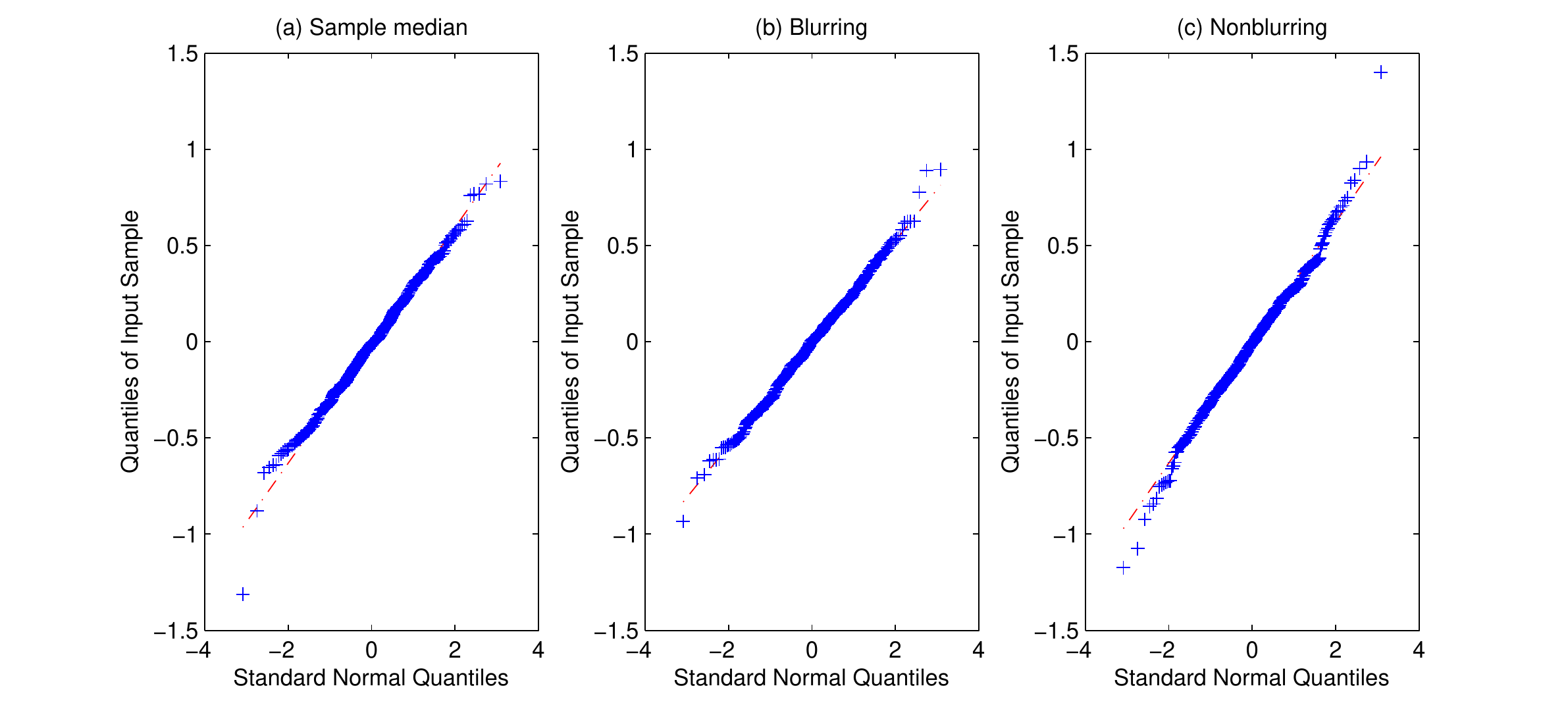}
\caption{QQ-plots. Sampling distributions are normal (top row), uniform (middle row)
and Student-$t$ (bottom row). Weight function is normal.}
\label{qqplot}
\end{figure}

\begin{figure}[htb]
$~~~~$\includegraphics[width=5.75in,height=2.5in]{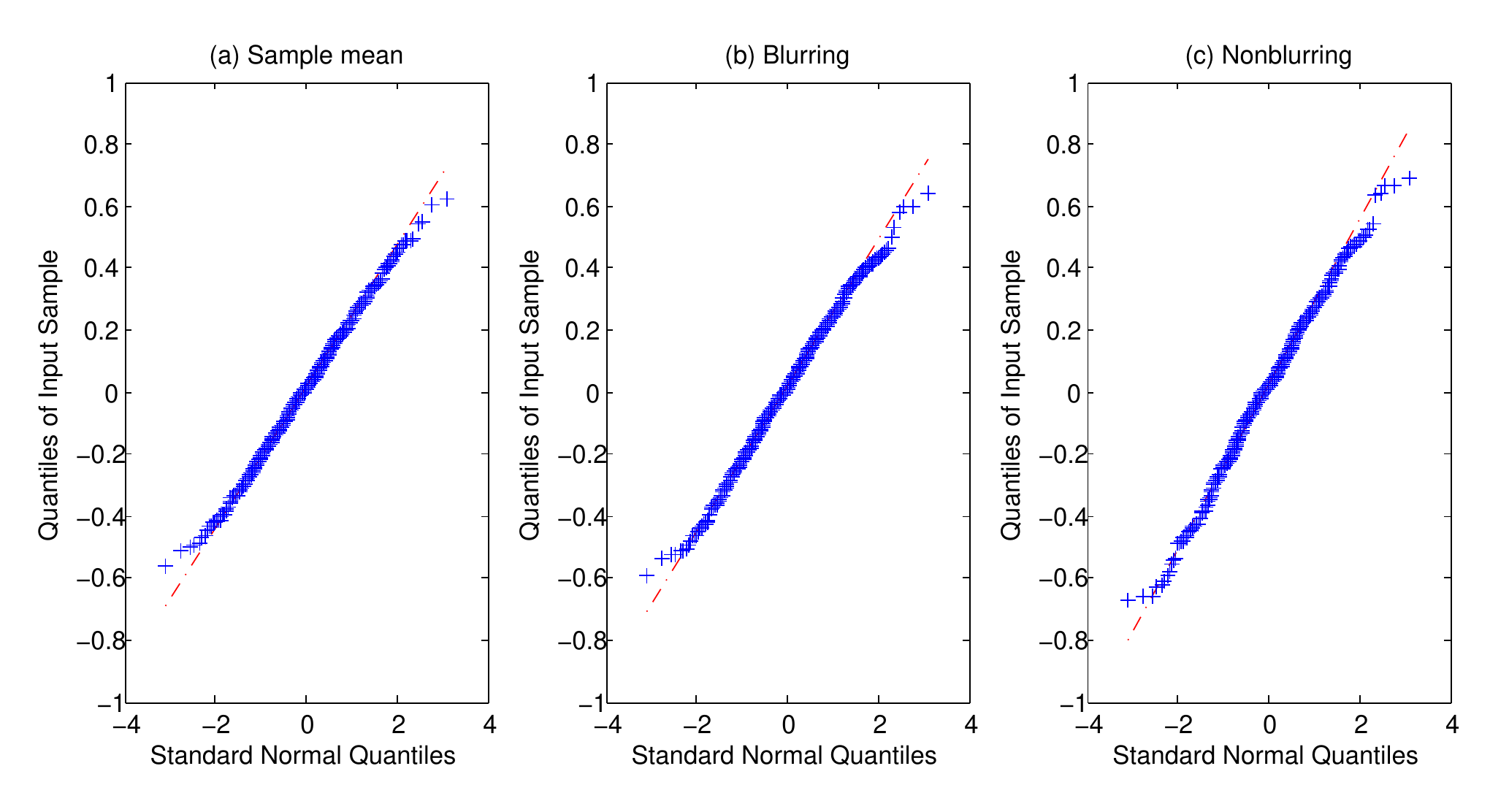}\\ $~$ \hspace{-0.75cm}
\includegraphics[width=6.5in,height=2.5in]{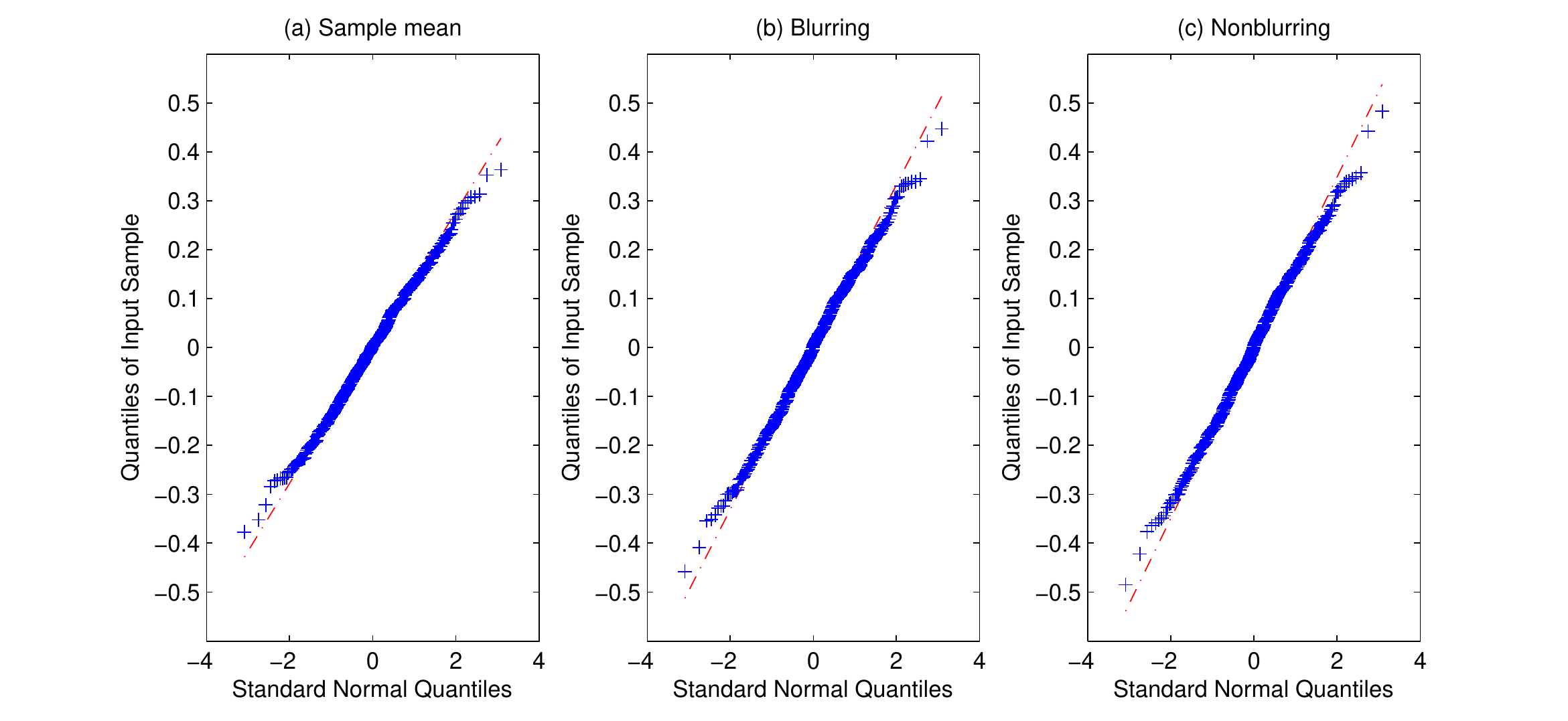}\\
$~~~~$\includegraphics[width=5.75in,height=2.5in]{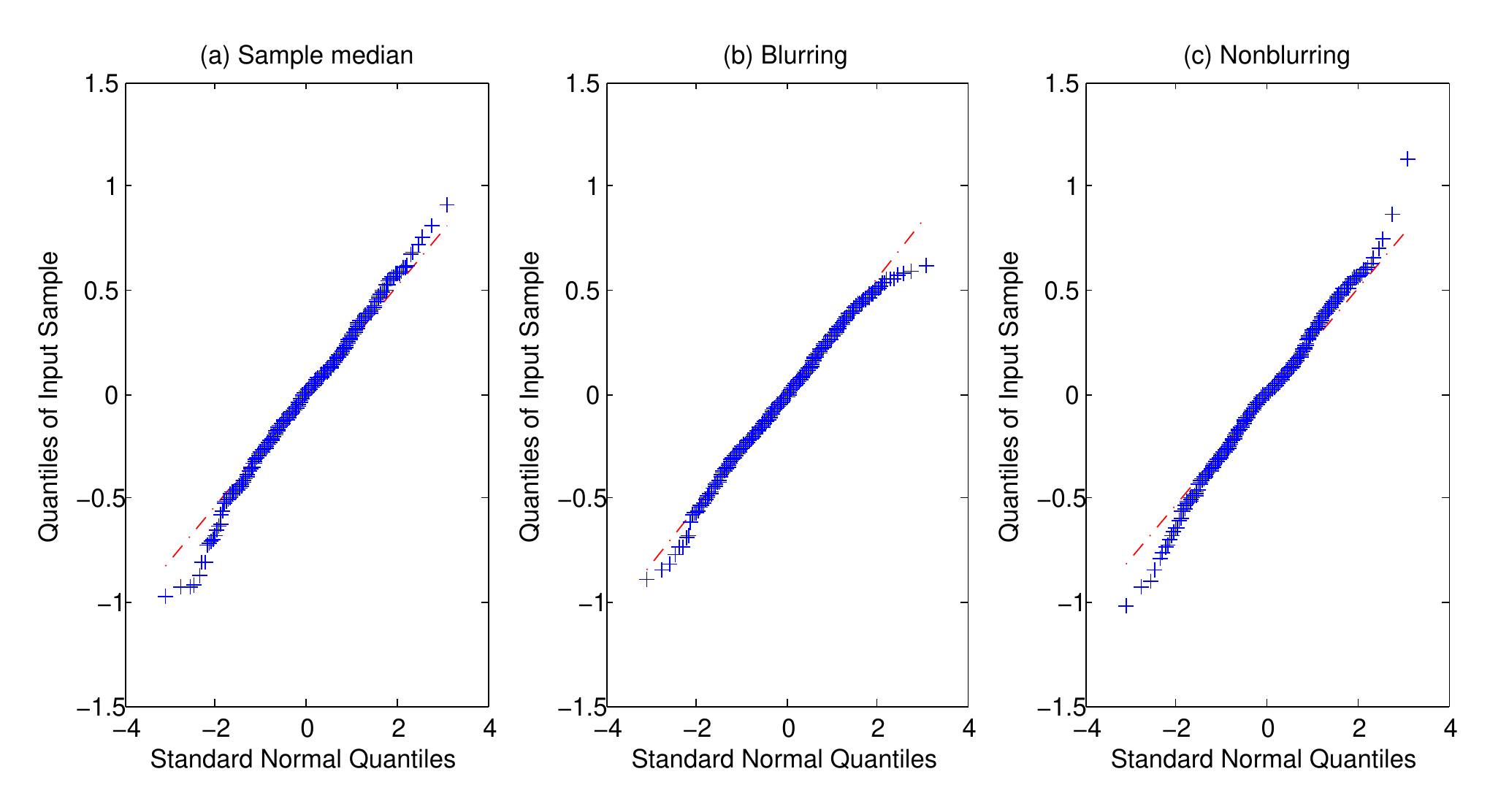}
\caption{QQ-plots. Sampling distributions are normal (top row), uniform (middle row)
and Student-$t$ (bottom row). Weight function is double exponential.}
\label{qqplot2}
\end{figure}

In this simulation study, we compare the asymptotic normality of
blurring and nonblurring processes. In theory all points should
converge to a common point \citep{Chen2} if the weight function has
unbounded support. However, in empirical data simulation (in
particular, the case of Student-$t$ distribution presented below)
some points far away from the main data cloud may fail to move to
the location where most points have converged to, due to the
precision in computer and the stopping criterion in the data
implementation. Therefore, we take the median value at the stopping
of the updating process. Precisely, let $\mu^\T_{n,m} :={\rm
median}\{x_{i,n,m}^\T,\, i=1,2,\dots,n\}$ and $\mu^\Tn_{n,m} :={\rm
median}\{x_{i,n,m}^\Tn,\, i=1,2,\dots,n\}$, where $\tau$ is the
number of iteration steps at convergence for the $m^{\rm th}$
replicate run for $m=1,2,\dots,M$. Here we take $M=500$. Data are
generated from $N(0,1)$, Uniform$(0,1)$, and Student-$t$ with 3
degrees of freedom, where the sample size $n$ is set to 400. Two
kinds of weight functions, normal and double exponential, are used.
In Figures~\ref{qqplot} and~\ref{qqplot2}, QQ-plots for
$\{\mu^\T_{n,m}:m=1,2,\dots, M\}$ and $\{\mu^\Tn_{n,m}:m=1,2,\dots,
M\}$ are presented for two kinds of weight functions. We also
include QQ-plots using sample means (for normal and uniform
distributions) and sample medians (for Student-$t$ distribution) as
a reference asymptotic behavior. For the Student-$t$, it requires a
much larger $n$ for the sample means to behave like a normal. Thus,
we used sample medians, which require less larger $n$. It can be
seen that both $\{\mu^\T_{n,m}\}_{m=1}^M$ and
$\{\mu^\Tn_{n,m}\}_{m=1}^M$ well follow a normal distribution.
Furthermore, the slopes in nonblurring QQ-plots are a bit steeper
than those in the blurring ones. It indicates that the convergence
point of nonblurring process tends to have a larger variance.

\subsection{Mean square error comparison}\label{compare_var}

In this simulation study, we compare mean square errors (MSE) of
$\{\mu^\T_{n,m}\}_{m=1}^M$ and $\{\mu^\Tn_{n,m}\}_{m=1}^M$, which
are obtained from blurring and nonblurring processes, respectively.
Three types of data distributions are used, $N(0,1)$, ${\rm
Uniform}(-1,1)$ and Student-$t$ with 3, 5 and 10 degrees of freedom.
Normal and double exponential weight functions are adopted.
Sample size is set to $n=10^2, 12^2, 14^2, \ldots, 30^2$. The mean
square errors are calculated based on $M=10,000$ multiple trials.

\begin{figure}[h]
\centering
\includegraphics[width=\textwidth,height=4in]{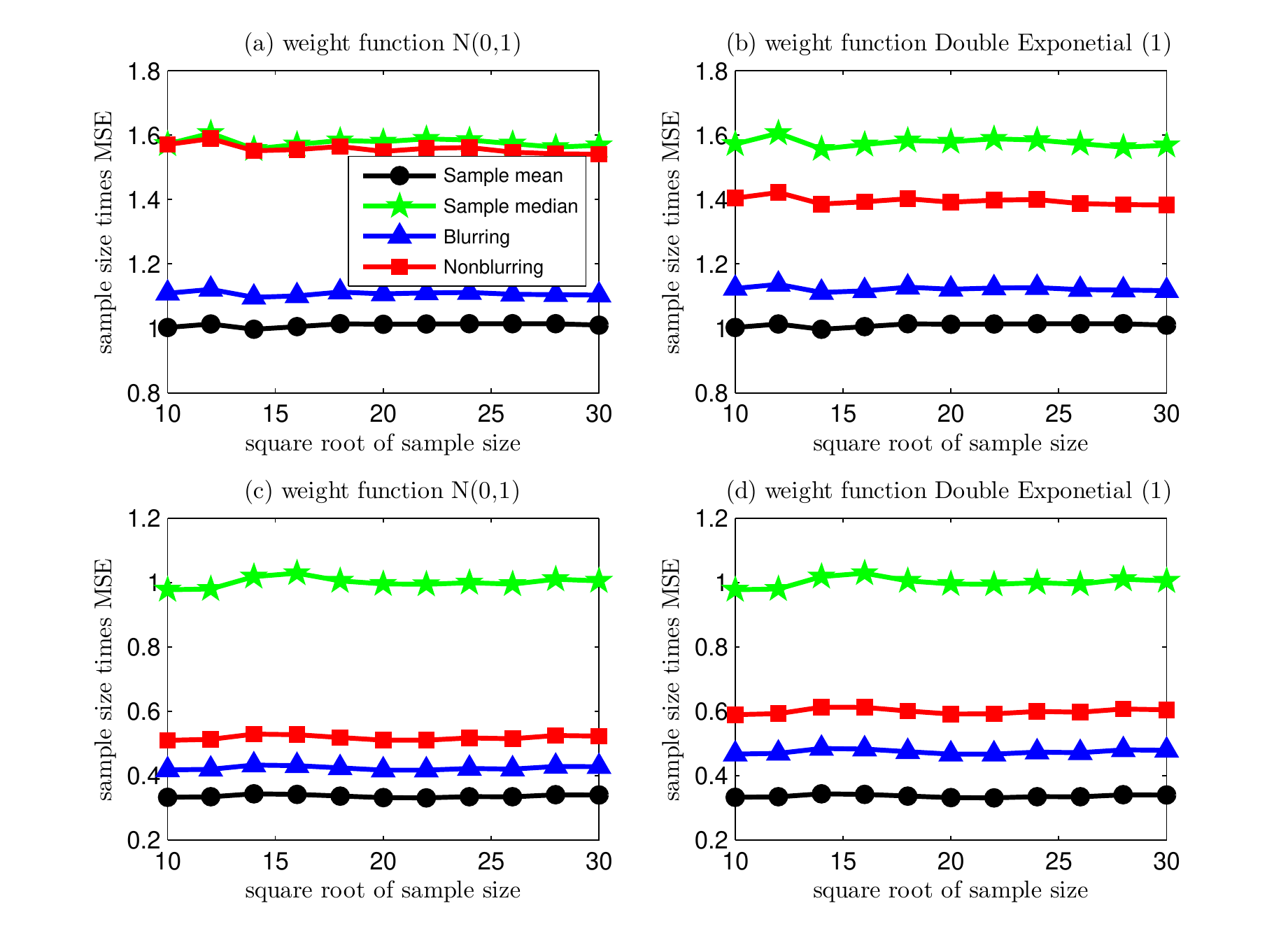}
\caption{MSE comparison.
(a) and (b): sampling distribution is normal $N(0,1)$; (c) and (d):
sampling distribution is uniform $(0,1)$}\label{mse1}
\end{figure}

Plots of $n \times $MSE against the square root of sample size
$\sqrt n$ are presented. We also include MSE of sample means and
sample medians for comparison. Sample mean is the UMVUE for
estimating the normal mean. From Figure~\ref{mse1} (a) and (b), we
can see that the sample mean has the smallest MSEs for data
generated from normal. Estimates by the nonblurring process have
larger MSEs compared to those by blurring process. We have
experimented with various parameter values of normal and double
exponential as the weight functions. All the results show that
\[\mbox{sample mean MSE $<$ blurring estimate MSE
$<$ nonblurring estimate MSE $<$ sample median MSE}.\] Similar phenomena can be observed
for data generated from the uniform $(0,1)$ distribution, which are
depicted in Figure~\ref{mse1} (c) and (d).

\begin{figure}[h]
\centering
\includegraphics[width=\textwidth,height=5in]{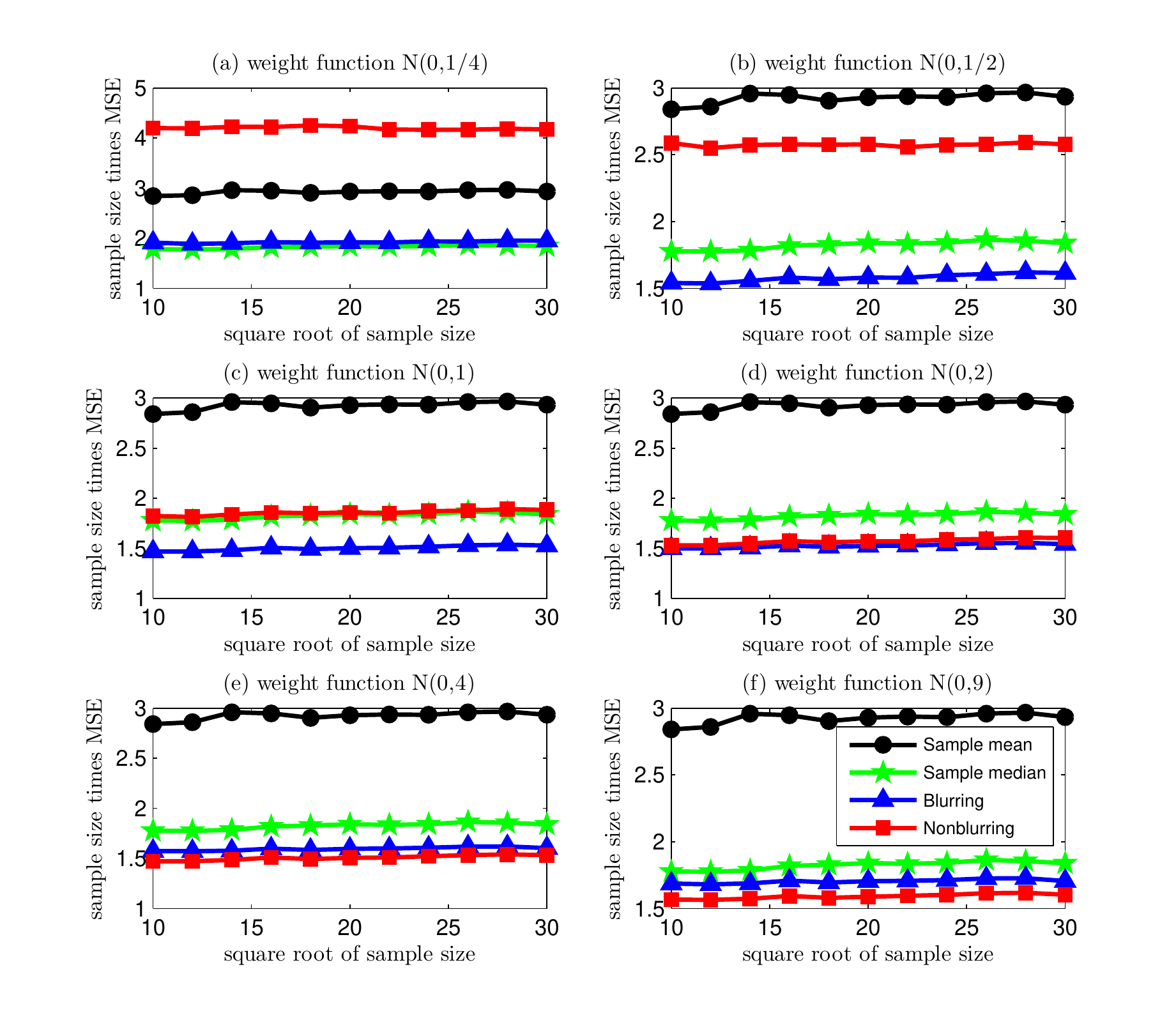}
\caption{MSE comparison for sampling distribution Student-$t$ with
3 degrees of freedom and normal weight function. }\label{mse2n}
\end{figure}

\begin{figure}[h]
\centering
\includegraphics[width=\textwidth,height=5in]{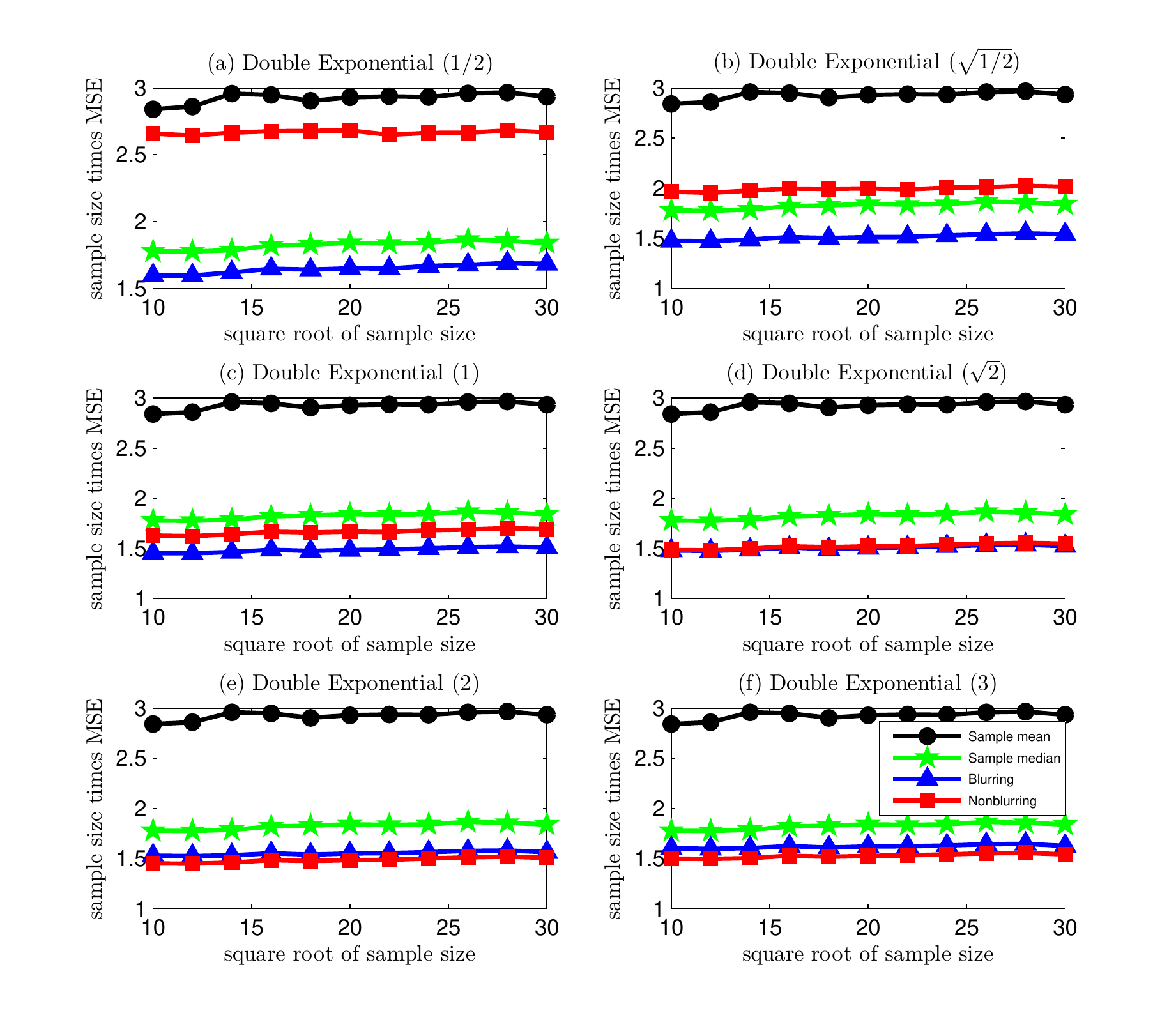}
\caption{MSE comparison for sampling distribution Student-$t$ with
3 degrees of freedom and double exponential weight function. }\label{mse2dex}
\end{figure}

For heavy-tailed distributions such as Student-$t$ distributions,
the performance of blurring estimate, nonblurring estimate, sample
mean and sample median for estimating location parameter depends on
the weight function. Figure~\ref{mse2n} presents the results by
normal weight functions with various parameter values, when the
sampling distribution is Student-$t$ with 3 degrees of freedom. As
expected, the sample median is better than the sample mean in
estimating the location parameter of Student-$t$ distribution. In
general, blurring estimates have smaller or competitive MSEs than
those by nonblurring estimates. They produce a bit larger MSEs than
the nonblurring ones only when an effectively flat weight function,
such as $N(0,\sigma^2)$ with $\sigma^2=4$ and $9$, is used. Similar
phenomena for double exponential weight functions on Student-$t$ can
be observed in Figure~\ref{mse2dex}.

\begin{figure}[h]
\centering
\includegraphics[width=\textwidth]{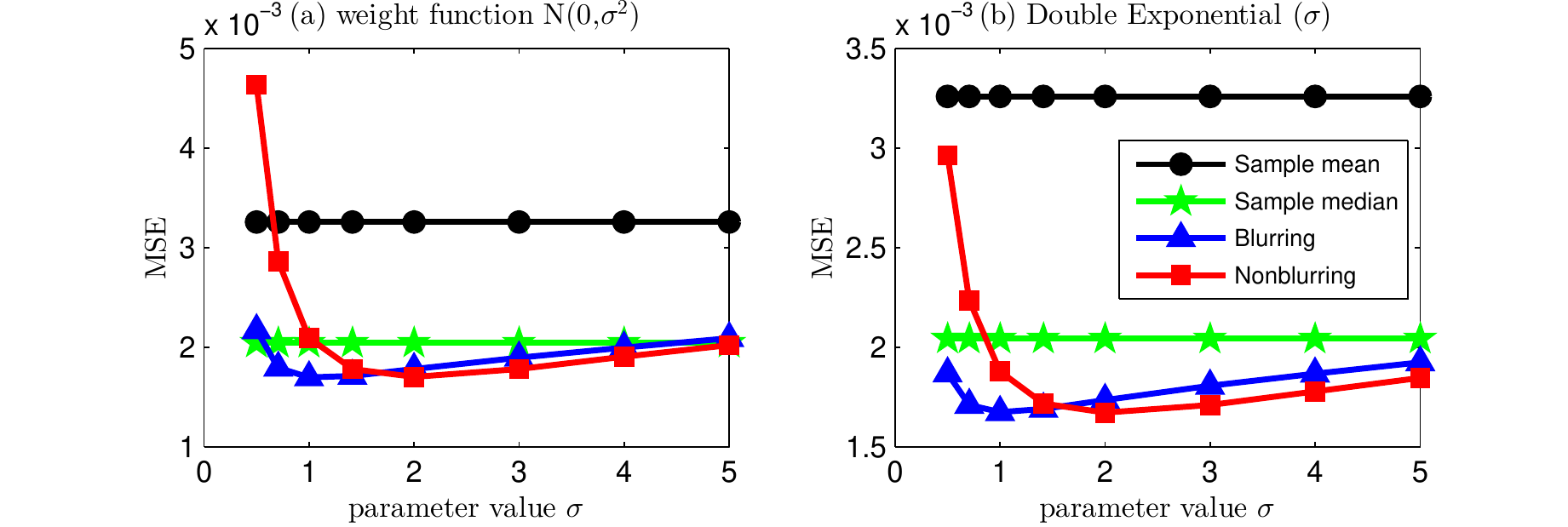}
\caption{Summary of MSE comparison for sampling distribution Student-$t$
with 3 degrees of freedom. }\label{mse2}
\end{figure}

\begin{figure}[h]
\centering
\includegraphics[width=\textwidth,height=5in]{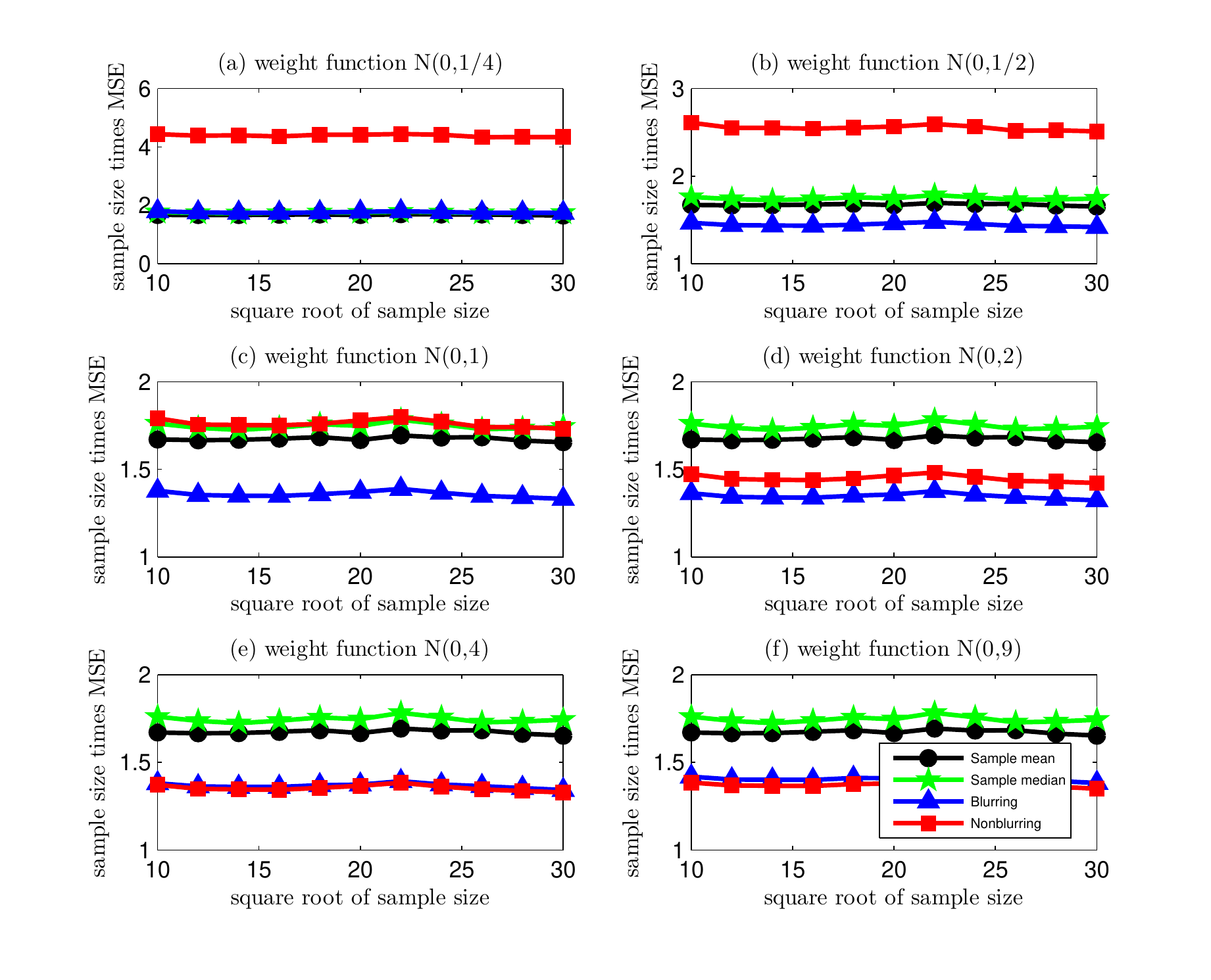}
\caption{MSE comparison for sampling distribution Student-$t$ with
5 degrees of freedom and normal weight function. }\label{mse5n}
\end{figure}

\begin{figure}[h]
\centering
\includegraphics[width=\textwidth]{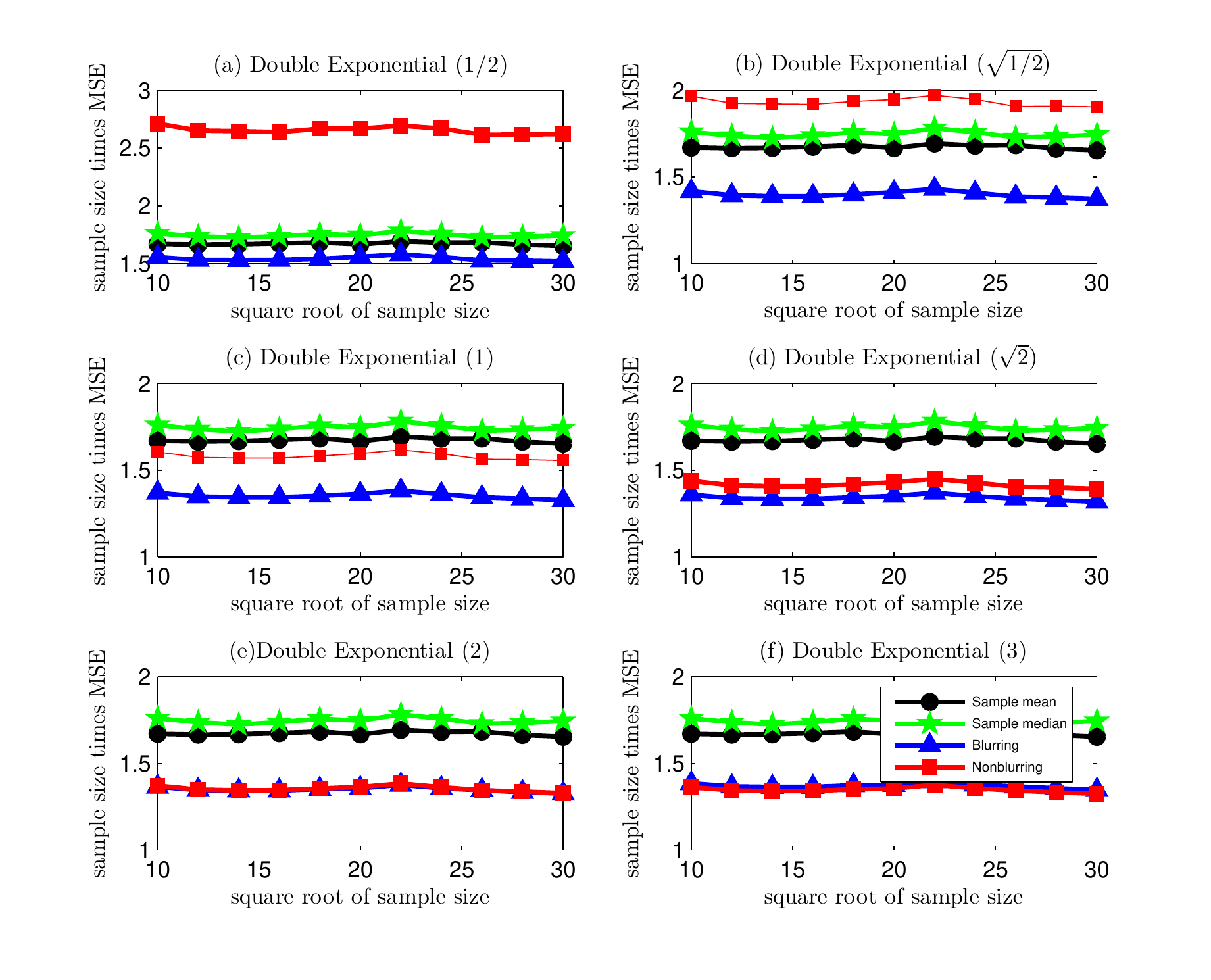}
\caption{MSE comparison for sampling distribution Student-$t$ with
5 degrees of freedom and double exponential weight function. }\label{mse5dex}
\end{figure}

\begin{figure}[h]
\centering
\includegraphics[width=\textwidth]{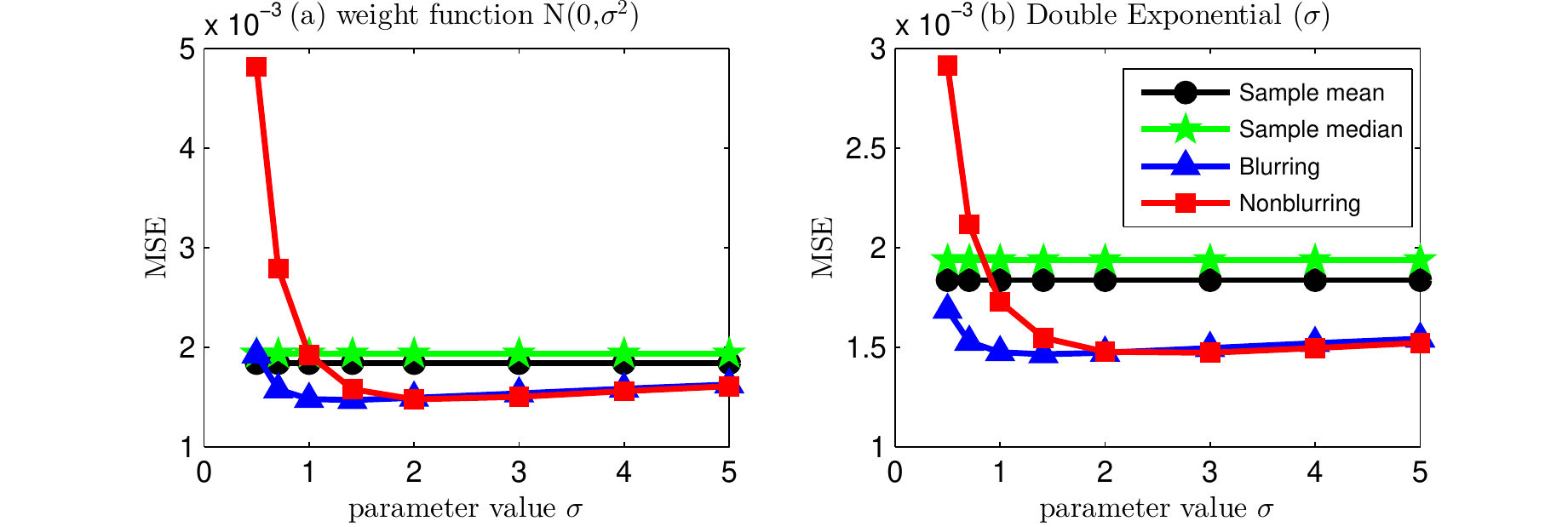}
\caption{Summary of MSE comparison for sampling distribution Student-$t$
with 5 degrees of freedom. }\label{mse5}
\end{figure}

\begin{figure}[h]
\centering
\includegraphics[width=\textwidth]{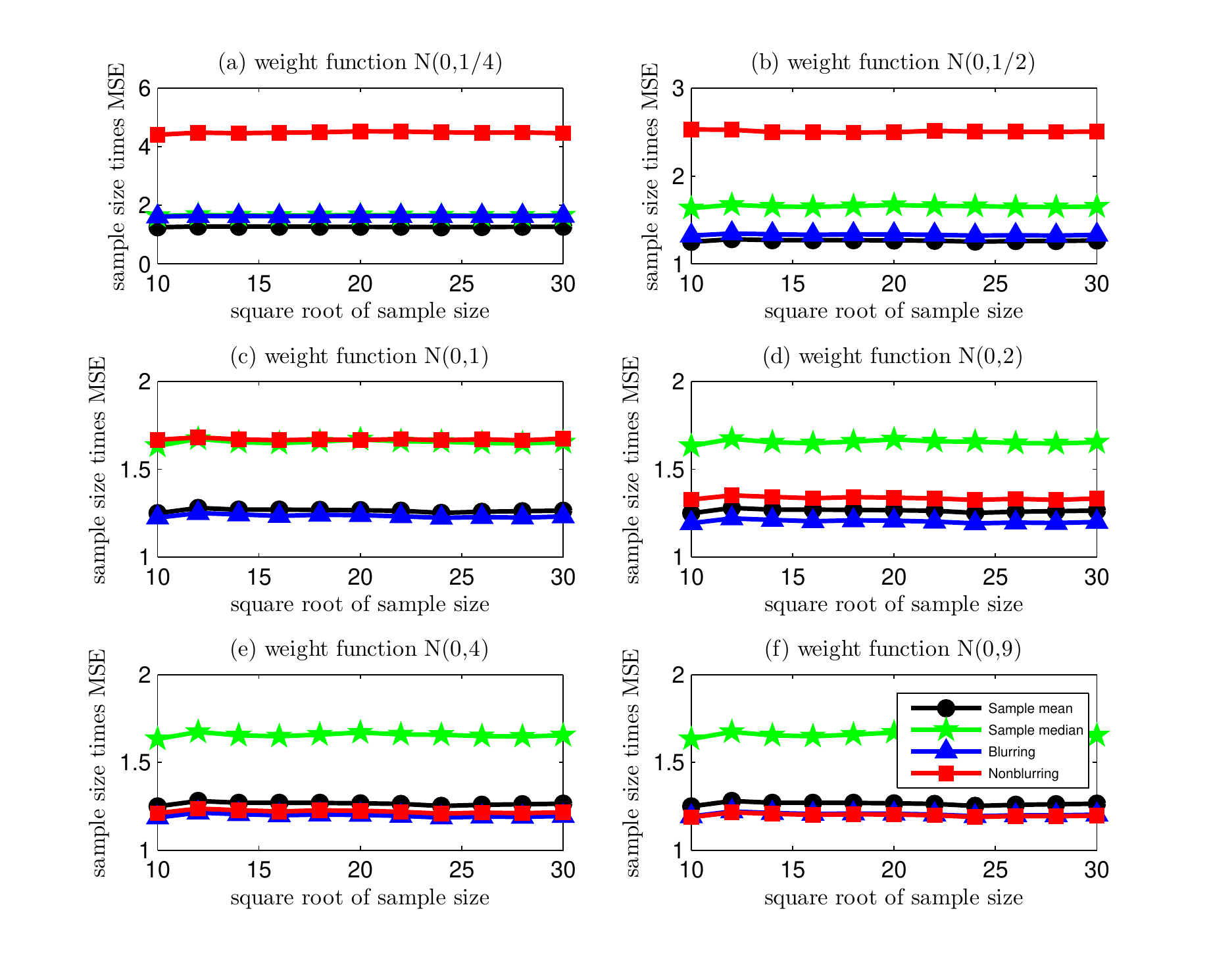}
\caption{MSE comparison for sampling distribution Student-$t$ with
10 degrees of freedom and normal weight function. }\label{mse10n}
\end{figure}

\begin{figure}[h]
\centering
\includegraphics[width=\textwidth]{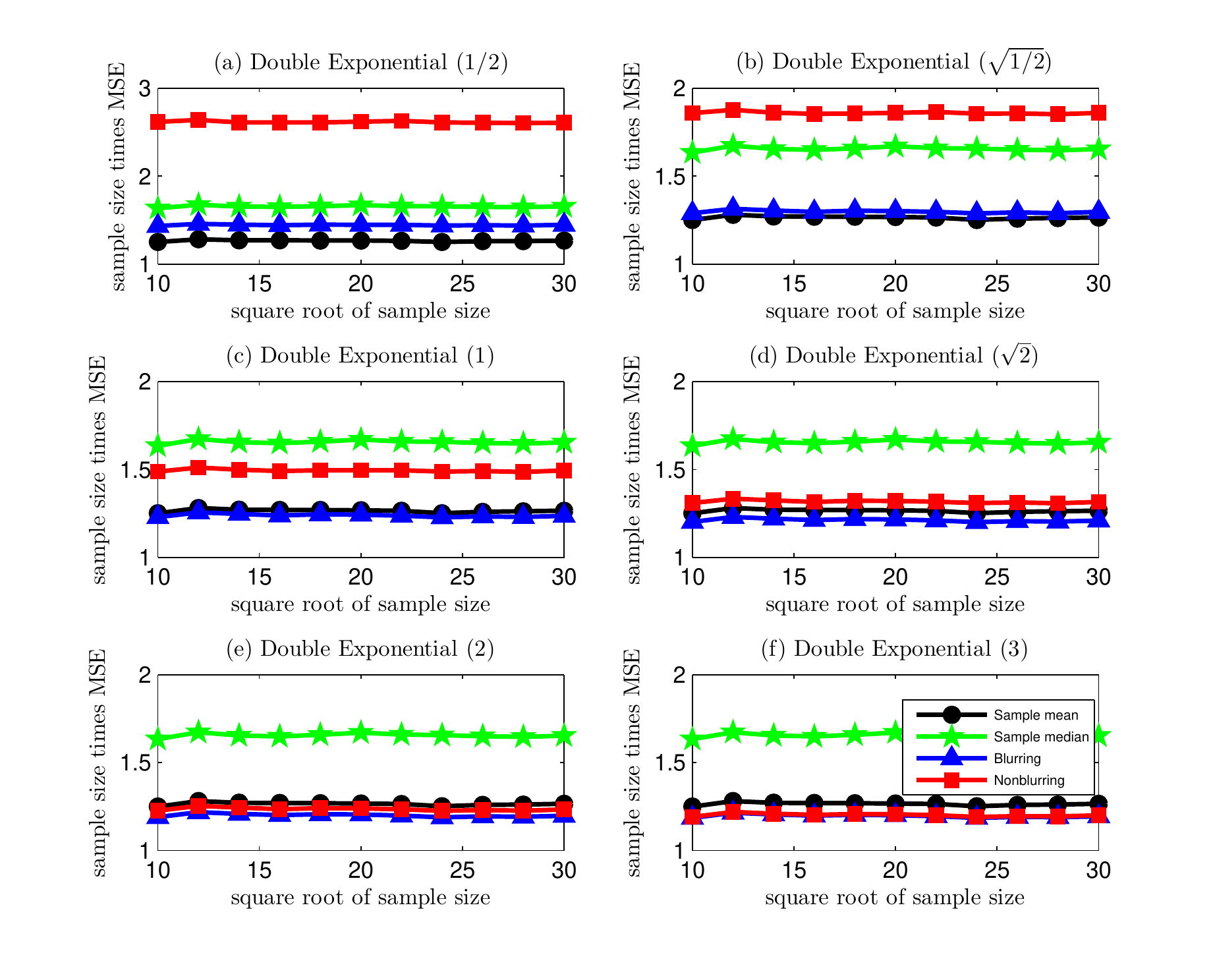}
\caption{MSE comparison for sampling distribution Student-$t$ with
10 degrees of freedom and double exponential weight function.
}\label{mse10dex}
\end{figure}

\begin{figure}[h]
\centering
\includegraphics[width=\textwidth]{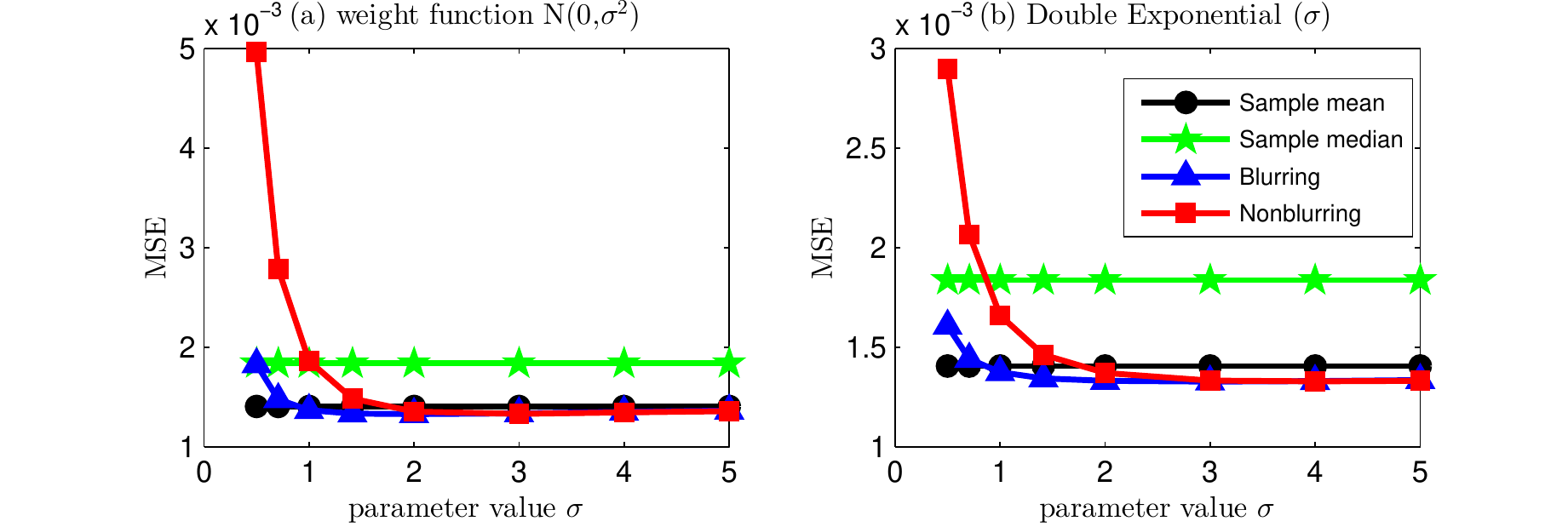}
\caption{Summary of MSE comparison for sampling distribution Student-$t$
with 10 degrees of freedom. }\label{mse10}
\end{figure}

To further compare MSEs of all estimates,  in Figure~\ref{mse2} we
summarize the MSE results with sample size $n=900$ on cases from
Figures~\ref{mse2n} and \ref{mse2dex}. From Figure~\ref{mse2}, the
best-performance weight functions for blurring estimates are
$N(0,1)$ and double exponential with mean parameter equal to 1. The
best-performance weight functions for nonblurring estimates are
$N(0,4)$ and double exponential with mean 2. All the MSE values are
very close in these 4 best-performance cases. While the nonblurring
estimation produces much larger MSEs for peaked weight functions,
the performance of blurring estimation is still competitive with
nonblurring for relatively flat weight function. Therefore, blurring
is viewed as a more robust estimator for Student-$t$. We have also
extended the experiment to Student-$t$ distributions with 5 and 10
degrees of freedom. The results are presented in
Figures~\ref{mse5n}-\ref{mse10}. They all support that the blurring
estimates are more robust in the sense of smaller or competitive
MSEs. The theoretic aspects behind these findings would be worthy of
further exploration.

\section{Conclusion}

In this article, we have established the weak convergence and
associated Central Limit Theorem for the blurring and nonblurring
processes. Convergence points from both types of processes can be
used for robust M-estimation of location.  In our simulation study,
it shows that the blurring type has a smaller mean square error than
the nonblurring type if a weight function is reasonably chosen. The
nonblurring type estimation is often adopted in robust statistics
literature. Our simulation results suggest that we shall consider
the blurring type algorithm as an alternative choice to the existent
nonblurring type algorithm for robust M-estimation.

\clearpage

\begin{appendix}
\section{Appendix}
Except for Lemma~\ref{lemma3} with a direct proof, we will show the
rest lemmas and Theorem~\ref{thm:main} by mathematical induction.
For each individual lemma or theorem, we will first show that it is
valid for $t=0$. Next by assuming the validity of
Lemmas~\ref{lemma:order}-\ref{lemma:sk_conv} as well as
Theorem~\ref{thm:main} for $t$ and the validity of preceding lemmas
for $t+1$, we will establish the claim of the current target lemma.
For instance, say, the current target lemma that we want to prove is
Lemma~\ref{lemma4}. We first show that Lemma~\ref{lemma4} is valid
for $t=0$. Next, by assuming that
Lemmas~\ref{lemma:order}-\ref{lemma:sk_conv} and
Theorem~\ref{thm:main} hold for $t$ and further assuming that
Lemmas~\ref{lemma:order}-\ref{lemma3} hold for $t+1$, we will
establish the claim of Lemma~\ref{lemma4} for $t+1$.

\subsection{Proof of Lemma \ref{lemma:order}}
\label{appendix:lemma1}
\begin{proof}
We will prove that for $0<a<b$,
\[
\eta^\1(a)=\frac {\int_{-\infty}^{\infty} x w(x-a)f(x) dx}
{\int_{-\infty}^{\infty}  w(x-a) f(x)dx} < \frac
{\int_{-\infty}^{\infty} x w(x-b)f(x) dx} {\int_{-\infty}^{\infty}
w(x-b) f(x) dx}=\eta^\1(b).
\]
for any symmetric probability density function $f$. The proof for
$\eta^\tplus(a) < \eta^\tplus(b)$ is identically the same. We first
consider the case that $\eta^\1(a) \leq \frac {a+b} 2$. From
\[
\frac {\int_{-\infty}^{\infty} x w(x-a)f(x) dx}
{\int_{-\infty}^{\infty}  w(x-a) f(x)dx}=\eta^\1(a),
\]
we have
\[
\int_{-\infty}^{\eta^\1(a)} (\eta^\1(a)-x) w(x-a)f(x) dx =\int_{\eta^\1(a)}^{\infty} (x-\eta^\1(a)) w(x-a)f(x) dx.
\]
Since $w$ is log-concave, for any $d>0$, $\log(w(x+d))-\log(w(x))$
is non-increasing for all $x$. Therefore
\[
\frac {w(x-a)} {w(x-b)} \geq \frac {w(\eta^\1(a)-a)} {w(\eta^\1(a)-b)}\equiv \gamma   \quad \quad \forall x< \eta^\1(a).
\]
and
\[
\frac {w(x-a)} {w(x-b)} \leq  \gamma   \quad \quad \forall x> \eta^\1(a).
\]
Then
\begin{equation}\label{eq:order1}
\int_{-\infty}^{\eta^\1(a)} (\eta^\1(a)-x) w(x-a)f(x) dx \geq \gamma \int_{-\infty}^{\eta^\1(a)} (\eta^\1(a)-x) w(x-b)f(x) dx.
\end{equation}
Similar, we have
\begin{equation}\label{eq:order2}
\int_{\eta^\1(a)}^{\infty} (x-\eta^\1(a)) w(x-a)f(x) dx < \gamma \int_{\eta^\1(a)}^{\infty} (x-\eta^\1(a)) w(x-b)f(x) dx.
\end{equation}
Note that $\gamma \geq 1$ since $\eta^\1(a) < \frac {a+b} 2$ . Since
$w$ is non-constant, there exists $x > \frac {a+b} 2$, such that
$\frac {w(x-a)} {w(x-b)} <1 \leq \gamma$. Therefore the inequality
in (\ref{eq:order2}) is strictly less. Combining (\ref{eq:order1})
and (\ref{eq:order2}), we have
\begin{eqnarray*}
&&\int_{-\infty}^{\infty} (x-\eta^\1(a)) w(x-b)f(x) dx\\
&=& \int_{-\infty}^{\eta^\1(a)} (x-\eta^\1(a)) w(x-b)f(x) dx
+ \int_{\eta^\1(a)}^{\infty} (x-\eta^\1(a)) w(x-b)f(x) dx\\
&>& \frac 1 \gamma \int_{-\infty}^{\eta^\1(a)} (x-\eta^\1(a)) w(x-a)f(x) dx
+ \frac 1 \gamma \int_{\eta^\1(a)}^{\infty} (x-\eta^\1(a)) w(x-a)f(x) dx\\
&=&\frac 1 \gamma \int_{-\infty}^{\infty} (x-\eta^\1(a)) w(x-a)f(x) dx =0
\end{eqnarray*}
Therefore
\[
 \eta^\1(b)-\eta^\1(a)=\frac
{\int_{-\infty}^{\infty} (x-\eta^\1(a)) w(x-b)f(x) dx} {\int_{-\infty}^{\infty}
w(x-b) f(x) dx}>0.
\]
For the case that $\eta^\1(a) > \frac {a+b}2$, the proof is almost
the same. Now
\[
\frac {w(x-b)} {w(x-a)} \geq \frac {w(\eta^\1(a)-b)}{w(\eta^\1(a)-a)}\equiv \theta  \quad \quad \forall x> \eta^\1(a).
\]
Then we have
\[
\int_{\eta^\1(a)}^{\infty} (x-\eta^\1(a)) w(x-b)f(x) dx \geq \theta \int_{\eta^\1(a)}^{\infty} (x-\eta^\1(a)) w(x-a)f(x) dx
\]
and
\[
\int_{-\infty}^{\eta^\1(a)} (\eta^\1(a)-x) w(x-b)f(x) dx < \theta \int_{-\infty}^{\eta^\1(a)} (\eta^\1(a)-x) w(x-a)f(x) dx.
\]
Combining both will again lead to $\eta^\1(b) -\eta^\1(a) >0$.

\end{proof}

\subsection{Proof of Lemma~\ref{lemma2}}\label{appendix:lemma2}
\begin{proof}
Similar to Lemma~\ref{lemma:order}, we will prove this lemma for
$\eta_n^\1$. The proof for $\eta_n^\tplus$ is identically the same.
By definition, we have
\begin{equation}\label{eq:g_n_splus}
\eta_n^\1(x) = \frac {\int y w(y-x) dF_n(y)} {\int w(y-x) dF_n(y)}.
\end{equation}
Since $w$ is finite integrable, $(y-x) w(y-x)$ is bounded. Then,
there exists a constant $M<\infty$ such that
\[
|y w(y-x)| \leq |(y-x) w(y-x)| + |x w(0)| \leq M +|x| w(0).
\]
From the weak convergence of $F_n$, we have
\begin{eqnarray}
\lim_{n \to \infty} \int y w(y-x) dF_n(y) &=& \int y w(y-x) dF(y) \quad \mbox{ a.s. }\label{eq:num}\\
\lim_{n \to \infty} \int w(y-x) dF_n(y) ~~ &=& \int w(y-x) dF(y) \quad \mbox{ a.s. }\label{eq:den}
\end{eqnarray}
Also note that, for a fixed $x$
\[
\int w(y-x) dF(y) > \int_0^x w(y-x) dF(y) > (F(x)-1/2)\times w(-x)
\] is bounded away from $0$. Thus, from (\ref{eq:g_n_splus}), (\ref{eq:num}) and (\ref{eq:den}),
\[
\lim_{n \to \infty} \eta_n^\1(x)=\eta^\1(x) \quad \mbox{ a.s. }
\]

\end{proof}

\subsection{Proof of Lemma~\ref{lemma3}}
\label{appendix:lemma3}
\begin{proof}
Let
\begin{eqnarray*}
\alpha_n&=& \int y w(y-x) dF_n^\t(y),\;\;
\alpha\, ~=\, \int y w(y-x) dF^\t(y), \\
\beta_n&=& \int  w(y-x) dF_n^\t(y), \;\;
~\, \beta~\, =\,~ \int  w(y-x) dF^\t(y).
\end{eqnarray*}
Then,
\begin{eqnarray*}
&&\sqrt{n}\left(\eta_n^\tplus(x)-\eta^\tplus(x)\right)
  = \sqrt{n}\left(\frac{\alpha_n}{\beta_n}-\frac{\alpha}{\beta}\right)\\
&=& \sqrt{n}\left(\frac{\alpha_n}{\beta}-\frac{\alpha}{\beta}+\frac{\alpha_n}{\beta_n}-\frac{\alpha_n}{\beta}\right)
  = \sqrt{n}\left(\frac{\alpha_n-\alpha}{\beta}-\frac{\alpha_n(\beta_n-\beta))}{\beta_n \beta}\right)\\
&=& \frac 1 \beta \int y w(y-x) d\sqrt{n}\left(F^\t_n(y)-F^\t(y)\right)\\
&& - \frac {\alpha_n}{\beta_n \beta}\int w(y-x) d\sqrt{n}\left(F^\t_n(y)-F^\t(y)\right)\\
&=& \frac {\int (y-\eta_n^\tplus(x)) w(y-x) d\sqrt{n}\left(F^\t_n(y)-F^\t(y)\right)}{\int  w(y-x) dF^\t(y)}.
\end{eqnarray*}
\end{proof}

\subsection{Proof of Lemma \ref{lemma4}}
\label{appendix:lemma4}

\begin{proof}
\begin{eqnarray}
&& Z_n^\tplus(x) =\sqrt{n}\left(F_n^\tplus(x)-F^\tplus(x)\right)\nonumber\\
&=&\sqrt{n}\left(F_n^\t(\xi_n^\tplus(x))
   - F^\t(\xi^\tplus(x) \right)\nonumber\\
&=&\sqrt{n}\left(F_n^\t(\xi_n^\tplus(x))
   - F_n^\t(\xi^\tplus(x) \right)
   +\sqrt{n}\left(F_n^\t(\xi^\tplus(x))
   - F^\t(\xi^\tplus(x) \right)\nonumber\\
&=& \sqrt{n}\left(F_n^\t(\xi_n^\tplus(x))- F_n^\t(\xi^\tplus(x) \right)
   + \int 1_{\{y\leq \xi^\tplus(x)\}} \, dZ_n^\t(x). \label{ineq:cdf2term}
\end{eqnarray}
For the first term, we need to to calculate
$\xi_n^\tplus(x)-\xi^\tplus(x)$. Since $\eta_n^\tplus(x)$ converges
to $\eta^\tplus(x)$ almost surely, the inverse function
$\xi_n^\tplus(x)$ also converges to $\xi^\tplus(x)$ almost surely.
Note that $x=\eta^\tplus(\xi^\tplus(x))$. Expanding
$\xi_n^\tplus(\cdot)$ at $\eta_n^\tplus(\xi^\tplus(x))$, we have
\begin{eqnarray*}
&&\xi_n^\tplus(x)-\xi^\tplus(x) \\
&=& \xi_n^\tplus\left(\eta_n^\tplus(\xi^\tplus(x))\right)
 - \frac {d  \xi_n^\tplus(u)}{d  u}\Big|_{u=\tau_n}
  \cdot \left(\eta_n^\tplus(\xi^\tplus(x))-x\right)-\xi^\tplus(x) \\
&=& - \frac {d  \xi_n^\tplus(u)}{d  u}\Big|_{u=\tau_n} \cdot
  \left(\eta_n^\tplus(\xi^\tplus(x))-\eta^\tplus(\xi^\tplus(x))\right)\\
&=& - \left(\frac {d  \xi^\tplus(u)}{d  u}\Big|_{u=x}+o_p(1)\right)
\cdot \left(\eta_n^\tplus(\xi^\tplus(x))-\eta^\tplus(\xi^\tplus(x))\right),
\end{eqnarray*}
where $\tau_n$ is some number between $x$ and
$\eta_n^\tplus(\xi^\tplus(x))$. Since $\eta_n^\tplus \to
\eta^\tplus$ a.s., $\eta_n^\tplus(\xi^\tplus(x)) \to x$ a.s.
Therefore, $\tau_n \to x$ a.s. Then the first term
in~(\ref{ineq:cdf2term}) can be calculated as follows.
\begin{eqnarray}
&&\sqrt{n}\left(F_n^\t(\xi_n^\tplus(x))   - F_n^\t(\xi^\tplus(x) \right)\nonumber\\[1ex]
&=& \frac { d  F_n^\t (\xi)} {d  \xi} \Big|_{\xi=c_n} \times \sqrt{n}\,
   \left(\xi_n^\tplus(x)-\xi^\tplus(x) \right)\nonumber\\[1ex]
&=& \left\{\frac {d  F^\t (\xi)} {d\xi} \Big|_{\xi=\xi^\tplus(x)} +o_p(1)\right\}\times
   \bigg\{-\left(\frac {d  \xi^\tplus(u)}{d  u}\Big|_{u=x}+o_p(1)\right) \nonumber\\
&&\times \sqrt{n}\left(\eta_n^\tplus(\xi^\tplus(x))
 -\eta^\tplus(\xi^\tplus(x))\right)\bigg\}\nonumber\\[1ex]
&=& -\left(\frac {d  F^\tplus (u)} {d  u} \Big|_{u=x} +o_p(1)\right)
 \times \sqrt{n}\left(\eta_n^\tplus(\xi^\tplus(x))
 -\eta^\tplus(\xi^\tplus(x))\right)\nonumber\\[1ex]
&=& -f^\tplus (x)\times \sqrt{n}\left(\eta_n^\tplus(\xi^\tplus(x))
 -\eta^\tplus(\xi^\tplus(x))\right) +o_p(1), \label{eta_n_in_lem3}
\end{eqnarray}
where $c_n$ is some number between $\xi_n^\tplus(x)$ and
$\xi^\tplus(x)$. Apply (\ref{eq:wc_g}) in Lemma~\ref{lemma3} to
(\ref{eta_n_in_lem3}), and then this lemma can be established.
\end{proof}

\subsection{Proof of Lemma \ref{lemma:sk_conv}}\label{appendix:lemma5}
\begin{proof}
Let $\Lambda$ denote the class of strictly increasing continuous
mappings from $[0,1]$ onto itself. For $\lambda \in \Lambda$, define
\[
\|\lambda\| := \sup_{0 \leq u \leq v \leq 1} \left| \log \left(\frac {\lambda(t)-\lambda(u)}{v-u}\right) \right|.
\]
The metric of $\cal D$ is defined by (see, e.g., Billingsley 1968)
\[
d(f_1,f_2)=\inf_\lambda\, \max \left\{ \|\lambda\|, \sup_{0\leq u \leq 1} \left|f_1(u)-f_2(\lambda (u))\right| \right\},~~
f_1,f_2 \in {\cal D}.
\]
The topology generated by this metric is called the Skorokhod
topology. For $d(f_1,f_2)<\delta$, there exists $\tilde{\lambda}$
such that
\[  \max \{ \|\tilde{\lambda}\|,
\sup_{0\leq u \leq 1} |f_1(u)-f_2(\tilde{\lambda} (u))| \}<\delta.\]
This implies that
\[
\|\tilde{\lambda}\| < \delta \qquad {\rm and}\qquad
\sup_{0\leq u \leq 1} |f_1(u)-f_2(\tilde{\lambda} (u))| <\delta.
\]
Then, we have
\begin{eqnarray*}
&&\left| ({\cal L}^\tplus f_1)(u)-({\cal L}^\tplus f_2)(\tilde{\lambda} (u)) \right|\\
&=& \left| \int_{v=0}^1 K^\tplus(F^\it(u),F^\it(v)) d \left( f_1(v)-f_2(\tilde{\lambda}(v))\right) \right| \\
&\leq& \delta \int_{v=0}^1 \left| K^\tplus(F^\it(u),F^\it(v)) \right| dv.
\end{eqnarray*}
Plugging in defining expression for $K^\tplus$, we have
\begin{eqnarray*}
&& \int_{v=0}^1 \left| K^\tplus(F^\it(u),F^\it(v)) \right| dv\\
&=&  \int_{y=-\infty}^\infty \left| K^\tplus(x,y) \right| d F^\t(y) \quad \quad \quad \quad \quad (x=F^\it(u), y=F^\it(v))\\[1ex]
&=& \int_{y=-\infty}^\infty \left|  -\frac{f^\tplus(x)\cdot(y-x) w(y-\xi^\tplus(x))}{ {\int  w\left(y-\xi^\tplus(x)\right) dF^\t(y)}}
+1_{\{y\leq \xi^\tplus(x)\}} \right| d F^\t(y)\\[1ex]
&\leq& \int_{y=-\infty}^\infty \left|  \frac{f^\tplus(x)\cdot |y+x|\cdot w(y-\xi^\tplus(x))}{ {\int  w\left(y-\xi^\tplus(x)\right) dF^\t(y)}}
 \right| d F^\t(y)+1.
\end{eqnarray*}
Note that
\[
\int_{y=-\infty}^\infty \left|  \frac{f^\tplus(x)\cdot x w(y-\xi^\tplus(x))}{ {\int  w\left(y-\xi^\tplus(x)\right) dF^\t(y)}}
 \right| d F^\t(y) = f^\tplus(x) |x|,
\]
and that
\begin{eqnarray*}
&& \int_{y=-\infty}^\infty \left|  \frac{y w(y-\xi^\tplus(x))} {\int  w\left(y-\xi^\tplus(x)\right) dF^\t(y)}
 \right| d F^\t(y)\\
&=& \int_{y=-|x|}^{|x|} \left|  \frac{y w(y-\xi^\tplus(x))}{ {\int  w\left(y-\xi^\tplus(x)\right) dF^\t(y)}}
 \right| d F^\t(y)
 +\int_{y>|x|} \left|  \frac{y w(y-\xi^\tplus(x))}{ {\int  w\left(y-\xi^\tplus(x)\right) dF^\t(y)}}\right| d F^\t(y)\\
&\leq& \int_{y=-\infty}^{\infty} \left|  \frac{x w(y-\xi^\tplus(x))}{ {\int  w\left(y-\xi^\tplus(x)\right) dF^\t(y)}}
 \right| d F^\t(y)
 +\int_{y>|x|} \left|  \frac{y w(y-\xi^\tplus(x))}{ {\int  w\left(y-\xi^\tplus(x)\right) dF^\t(y)}} \right| d F^\t(y)\\
 &=& |x| + R^\t(x),
\end{eqnarray*}
where
\[
R^\t(x) := \int_{y>|x|} \left|  \frac{y w(y-\xi^\tplus(x))}{ {\int  w\left(y-\xi^\tplus(x)\right) dF^\t(y)}}
 \right| d F^\t(y) \to 0 \quad \mbox{ as } |x| \to \infty.
\]
Therefore, for $d(f_1,f_2)<\delta$ and $f_1,f_2\in{\cal D}[0,1]$, we
have
\begin{eqnarray*}
d({\cal L} f_1,{\cal L}f_2) &\leq& \max \left\{ \|\tilde{\lambda}\|, \sup_{0\leq u \leq 1} | ({\cal L} f_1)(u)-({\cal L} f_2)(\tilde{\lambda} (u))|\right\} \\
&<& \delta \left\{f^\tplus(x)(2|x|+R^\t(x))+1\right\},~~ {\rm where~} x=F^\it(u).
\end{eqnarray*}
Hence, ${\cal L}$ is continuous.
\end{proof}

\subsection{Proof of Theorem~\ref{thm:main} }\label{appendix:thm1}
\begin{proof} We will prove this theorem by mathematical induction.
For $t=0$, statements (i) and (ii) are well-known results as almost
sure convergence of empirical CDF and Donsker's Theorem,
respectively. (See, e.g., Dudley, 1999.) Assume statements (i) and
(ii) hold for $t=s$. Then, statement~(i) with $t=s+1$ is an
immediate result of Lemma~\ref{lemma2}. It is now left to show
statement~(ii) with $t=s+1$ to complete the proof by mathematical
induction.

Recall that $B_n^\s(u)= Z_n^\s(F^\is(u))$. By
Lemma~\ref{lemma:sk_conv},
\[
B^\splus(u)= \int_{v=0}^1 K^\splus(F^\is(u),F^\is(v)) dB^\s(v).
\]
By assumption $B^\s(u)= \int H^\s(u,z) dB^{(0)}(z)$, then
\[
d B^\s(u) = \int  \frac {\partial H^\s(u,z)}{\partial u}\, dB^{(0)}(z)\,du.
\]
Therefore,
\begin{eqnarray*}
B^\splus(u)&=& \int_{v=0}^1 K^\splus(F^\is(u),F^\is(v)) dB^\s(v)\\
&=&\int_{v=0}^1 \int_{z=0}^1 K^\splus(F^\is(u),F^\is(v)) \frac {\partial H^\s(v,z) }{\partial v}\, dB^{(0)}(z)\, dv \\
&=&\int_{z=0}^1 \int_{v=0}^1  K^\splus(F^\is(u),F^\is(v)) \frac {\partial H^\s(v,z) }{\partial v}\, dv\, dB^{(0)}(z) \\
&=&\int_{z=0}^1 H^\splus(u,z) dB^{(0)}(z),
\end{eqnarray*}
where $H^\1 (u,z):= K^\1(F^{-1}(u),F^{-1}(z))$ and
\[
H^\splus(u,z):= \int_{v=0}^1 K^\splus(F^\is(u),F^\is(v)) \frac {\partial H^\s(v,z) }{\partial v}\, dv.
\]
By mathematical induction, we have shown the almost sure convergence
of $F_n^\t$ and the weak convergence of
$\sqrt{n}\left(F_n^{(t)}(x)-F^{(t)}(x)\right)$.
\end{proof}

\subsection{Proof of Theorem~\ref{thm:main2} }\label{appendix:thm2}
\begin{proof}
We will use similar mathematical induction arguments to prove the
weak convergence for the nonblurring case. For $t=1$, nonblurring
and blurring process are identically the same. Therefore the
statements hold for $t=1$. Assume that they hold for $t=s$, we will
prove that they hold for $t=s+1$. Since
\[
\eta^\sn(x)=\frac {\int y \cdot w(y-x) dF(y)} {\int w(y-x) dF(y)},
\]
$\eta^\sn(x)$'s are the same for all $t$, i.e. $\eta^{[1]}(x)
=\eta^{[2]}(x)=\cdots$. With the same arguments for blurring
process, we have
\[
\lim_{n \to \infty} \eta_n^\snplus(x) =\eta^\snplus(x),
\]
and
\[
\sqrt{n}\left(\eta_n^\snplus(x)-\eta^\snplus(x)\right)=\frac {\int (y-\eta_n^\snplus(x)) w(y-x) d\sqrt{n}\left(F_n(y)-F(y)\right)}{\int  w(y-x) dF(y)}.
\]
Then by similar arguments,
\begin{eqnarray}
&&\sqrt{n}\left(F_n^\snplus(x)-F^\snplus(x)\right)\nonumber\\[0.5ex]
&=&\sqrt{n}\left(F_n^\sn(\xi_n^\snplus(x))
   - F_n^\sn(\xi^\snplus(x) \right)+\sqrt{n}\left(F_n^\sn(\xi^\snplus(x))
   - F^\sn(\xi^\snplus(x) \right)\nonumber\\[1ex]
&=& \left(\frac {d  F^\sn (\xi)} {d\xi} \Big|_{\xi=\xi^\snplus(x)} +o_p(1)\right)\cdot
 \left(-\frac {d  \xi^\snplus(u)}{d  u}\Big|_{u=x}+o_p(1)\right) \nonumber\\
&&\times \sqrt{n}\left(\eta_n^\snplus(\xi^\snplus(x))-\eta^\snplus(\xi^\snplus(x))\right)\nonumber\\
&&+\sqrt{n}\left(F_n^\sn(\xi^\snplus(x))
   - F^\sn(\xi^\snplus(x) \right)\nonumber\\[1ex]
&=&\int_y (K^\snplus(x,y)+o_p(1)) \, d\sqrt{n}\left(F_n(y)-F(y)\right)\nonumber \\
&&+\sqrt{n}\left(F_n^\sn(\xi^\snplus(x))
   - F^\sn(\xi^\snplus(x) \right),\label{eq:nb_part}
\end{eqnarray}
where
\[
K^\snplus(x,y)= -\frac{f^\snplus(x)\cdot(y-x) w(y-\eta^\isnplus(x))}{ {\int  w(y-\eta^\isnplus(x)) dF(y)}}.
\]
For $t=1$,
\begin{eqnarray*}
&&\sqrt{n}\left(F_n^{[1]}(x)-F^{[1]}(x)\right)\\
&=&\int_y (K^{[1]}(x,y)+o_p(1)) \, d\sqrt{n}\left(F_n(y)-F(y)\right)
  + \sqrt{n} \left(F_n^{[0]}(\xi^{[1]}(x))
   - F^{[0]}(\xi^{[1]}(x) \right)\\
&=& \int_y (K^{[1]}(x,y)+1_{\{y \leq \xi^{[1]}(x)\}}+ o_p(1)) \, d\sqrt{n}\left(F_n(y)-F(y)\right)\\
&=& \int_y (K^{[1]}(x,y)+1_{\{\eta^{[1]}(y) \leq x\}} + o_p(1)) \, d\sqrt{n}\left(F_n(y)-F(y)\right).
\end{eqnarray*}
Let $H^{[1]}(u,v)=K^{[1]}(F^{-1}(u),F^{-1}(v)) +
1_{\{\eta^{[1]}(F^{-1}(v)) \leq F^{-1}(u)\}}$. By mathematical
induction assumption of the statement at $t=s$, we have
\begin{eqnarray*}
&&\sqrt{n}\left(F_n^\sn(\xi^\snplus(x))-F^\sn(\xi^\snplus(y))\right) \\
&=& \int_y \left\{H^\sn \left(F(\xi^\snplus(x)),F(\xi^\snplus(y))\right)+o_p(1)\right\}
 d\sqrt{n}\left(F_n(y)-F(y)\right),
\end{eqnarray*}
where $H^\tn$, $t=1,2,\dots$, is defined iteratively
via~(\ref{eq:h_iter}). Take this into (\ref{eq:nb_part}), it becomes
\begin{eqnarray*}
&&\sqrt{n}\left(F_n^\snplus(x)-F^\snplus(x)\right)\\
& =& \int_y \left\{ K^\snplus (x,y) +   H^\sn \left(F(\xi^\snplus(x)),F(\xi^\snplus(y))\right)
 +o_p(1)\right\} d\sqrt{n}\left(F_n(y)-F(y)\right)\\
&=&  \int_y (H^\snplus(F(x),F(y))+o_p(1)) d\sqrt{n}\left(F_n(y)-F(y)\right),
\end{eqnarray*}
where
\begin{equation} \label{eq:h_iter}
H^\snplus(u,v) := K^\snplus(F^{-1}(u)),F^{-1}(v))
 +   H^\sn \left(F(\xi^\snplus(F^{-1}(u))),F(\xi^\snplus(F^{-1}(v)))\right).
\end{equation}
Then by similar arguments as in the blurring case, we have
\[
\sqrt{n}\left(F_n^\snplus(x)-F^\snplus(x)\right) \rightsquigarrow B^\snplus(F(x)),
\]
where
\[
B^\snplus(u)=\int_{v=0}^1 H^\snplus(u,v) dB^{[0]}(v).
\]
\end{proof}

\subsection{Proof of Theorem~\ref{thm:clt_bl}}\label{appendix:thm3}
\begin{proof}
From
\[
\rho^\t(x_{i,n}^\t)=\frac 1n \sum_{j=1}^n w(x_{j,n}^\t-x_{i,n}^\t)+o_p(1),
\]
then
\begin{eqnarray*}
&&\frac{\sum_{j=1}^n \frac 1n w(x_{j,n}^\t-x_{i,n}^\t) x_{j,n}^\t}{\sum_{j=1}^n \frac 1n w(x_{j,n}^\t-x_{i,n}^\t)}-
\frac{\sum_{j=1}^n \frac 1n w(x_{j,n}^\t-x_{i,n}^\t) x_{j,n}^\t}{\rho^\t(x_{i,n}^\t)}\\
&=&\frac{\sum_{j=1}^n \frac 1n w(x_{j,n}^\t-x_{i,n}^\t) x_{j,n}^\t}{\sum_{j=1}^n \frac 1n w(x_{j,n}^\t-x_{i,n}^\t)}\cdot  \frac{\rho^\t(x_{i,n}^\t)-\sum_{j=1}^n \frac 1n w(x_{j,n}^\t-x_{i,n}^\t)}{\rho^\t(x_{i,n}^\t)}\\
&=&o_p(1),
\end{eqnarray*}
since $\frac{\sum_{j=1}^n  w(x_{j,n}^\t-x_{i,n}^\t)
x_{j,n}^\t}{\sum_{j=1}^n  w(x_{j,n}^\t-x_{i,n}^\t)}$ is bounded.
Therefore
\begin{eqnarray}
S_n^\tplus &=&\frac1{\sqrt n}\sum_{i=1}^n x_{i,n}^\tplus
 =\frac1{\sqrt n}\sum_{i=1}^n
\frac{\sum_{j=1}^n w(x_{j,n}^\t-x_{i,n}^\t) x_{j,n}^\t}{\sum_{j=1}^n w(x_{j,n}^\t-x_{i,n}^\t)}\nonumber\\
&=&\frac1{n\sqrt n}\sum_{i,j=1}^n \frac {w(x_{j,n}^\t-x_{i,n}^\t)x_{j,n}^\t}{\rho^\t(x_{i,n}^\t)} +o_p(1)\nonumber\\
&=&\sqrt{n}\int\int \frac{w(x-y)x}{\rho^\t(y)} dF_n^\t(x) dF_n^\t(y) +o_p(1). \label{eq:Sn1}
\end{eqnarray}
Next, we have
\begin{eqnarray}
&&\sqrt{n}\left( \int\int \frac{w(x-y)x}{\rho^\t(y)} dF_n^\t(x)
dF_n^\t(y)- \int\int \frac{w(x-y)x}{\rho^\t(y)} dF^\t(x) dF^\t(y)
\right)\nonumber\\
&=&\sqrt{n}\left( \int\int \frac{w(x-y)x}{\rho^\t(y)} dF_n^\t(x)
dF_n^\t(y)- \int\int \frac{w(x-y)x}{\rho^\t(y)} dF_n^\t(x) dF^\t(y)
\right)\nonumber\\
&&+\sqrt{n}\left( \int\int \frac{w(x-y)x}{\rho^\t(y)} dF_n^\t(x)
dF^\t(y)- \int\int \frac{w(x-y)x}{\rho^\t(y)} dF^\t(x) dF^\t(y)
\right). \label{eq:Sn2}
\end{eqnarray}
By Fubini theorem, since the integral of the absolute value is
finite, we can change the order.
\begin{eqnarray*}
(\ref{eq:Sn2}) &=&\sqrt{n}\left( \int\int \frac{w(x-y)x}{\rho^\t(y)} dF_n^\t(x)
dF_n^\t(y)- \int\int \frac{w(x-y)x}{\rho^\t(y)} dF_n^\t(x) dF^\t(y)
\right)\nonumber\\
&&+\sqrt{n}\left( \int\int \frac{w(x-y)x}{\rho^\t(y)}
dF^\t(y)dF_n^\t(x)- \int\int \frac{w(x-y)x}{\rho^\t(y)} dF^\t(y)
dF^\t(x) \right)\nonumber\\
&=&\int\int \frac{w(x-y)x}{\rho^\t(y)} dF_n^\t(x)
d\sqrt{n}(F_n^\t(y)-F^\t(y))\\
&&+\int\int \frac{w(x-y)x}{\rho^\t(y)} dF^\t(y)d\sqrt{n}(F_n^\t(x)-
F^\t(x))\\
&\to& \int\left(\int (\frac{w(x-y)x}{\rho^\t(y)} +
\frac{w(y-x)y}{\rho^\t(x)}) dF^\t(x)\right)dB^\t(F^\t(y)).
\end{eqnarray*}
Since
\[
\int\int \frac{w(x-y)x}{\rho^\t(y)} dF^\t(x) dF^\t(y)=0,
\]
then
\begin{eqnarray*}
 S_n^\tplus
&=&S_n^\tplus - \sqrt{n}\int\int \frac{w(x-y)x}{\rho^\t(y)} dF^\t(x) dF^\t(y)\\
&\to& \int\int \left(\frac{w(x-y)x}{\rho^\t(y)} +
\frac{w(y-x)y}{\rho^\t(x)}\right) dF^\t(x) dB^\t(F^\t(y)).
\end{eqnarray*}
\end{proof}
\end{appendix}


\begin{thebibliography}{99}
\expandafter\ifx\csname
natexlab\endcsname\relax\def\natexlab#1{#1}\fi

\bibitem[Billingsley, 1968]{Billingsley}
Billingsley, P. (1968). {\em Convergence of Probability Measures},
New York, Wiley.

\bibitem[Basu et~al., 1998]{Basu}
Basu, A., Harris, I.~R., Hjort, N.~L., and Jones, M.~C. (1998).
\newblock Robust and efficient estimation by minimising a density power
  divergence.
\newblock {\em Biometrika}, 85(3):549--559.

\bibitem[Chen, 2015]{Chen2}
Chen, T.-L. (2015).
\newblock On the convergence and consistency of the blurring mean-shift
  process. {\em Annals of the Institute of Statistical
  Mathematics}, 67(1):157-176.

\item
Chen, T.L., Hsieh, D.N., Hung, H., Tu, I.P., Wu, P.S., Wu, Y.M.,
Chang, W. and Huang, S.Y. (2014). $\gamma$-SUP: a clustering
algorithm for cryo-electron microscopy images of asymmetric
particles. {\em Annals of Applied Statistics}, 8(1):259--285.

\bibitem[Chen and Shiu, 2007]{Chen}
Chen, T.-L. and Shiu, S.-Y. (2007).
\newblock A clustering algorithm by self-updating process.
\newblock {\em JSM Proceedings}, Statistical Computing Section, Salt Lake City,
  Utah; American Statistical Association, pp:2034--2038.

\bibitem[Cheng, 1995]{ChengY}
Cheng, Y. (1995).
\newblock Mean shift, mode seeking, and clustering.
\newblock {\em IEEE Transactions on Pattern Analysis and Machine Intelligence},
  17:790--799.

\bibitem[Comaniciu and Meer, 2002]{Comaniciu}
Comaniciu, D. and Meer, P. (2002). Mean shift: a robust approach
toward feature space analysis. {\em IEEE Transactions on Pattern
Analysis and Machine Intelligence}, 24(5):603-619.

\item
Dudley, R.M. (1999). {\em Uniform Central Limit Theorems}. Cambridge
University Press.

\bibitem[Field and Smith, 1994]{Field}
Field, C. and Smith, B. (1994).
\newblock Robust estimation: a weighted maximum likelihood approach.
\newblock {\em International Statistical Review}, 62(3):405--424.

\bibitem[Fujisawa and Eguchi, 2008]{Fujisawa}
Fujisawa, H. and Eguchi, S. (2008).
\newblock Robust parameter estimation with a small bias against heavy
  contamination.
\newblock {\em Journal of Multivariate Analysis}, 99:2053--2081.

\bibitem[Fukunaga and Hostetler, 1975]{Fukunaga}
Fukunaga, K. and Hosterler, L. D. (1975). The estimation of the
gradient of a density function, with applications in pattern
recognition. {\em IEEE Trans. Inform. Theory}, 21(1):32--40.

\bibitem[Ghassabeh, 2015]{Ghassabeh}
Ghassabeh, Y.A. (2015). A sufficient condition for the convergence
of the mean shift algorithm with Gaussian kernel. {\em Journal of
Multivariate Analysis}, 135:1-10.

\item
Hampel, F.R., Ronchetti, E.M., Rousseeuw, P.J. and Stahel, W.A.
(1986). {\em Robust Statistics: The Approach Based on Influence
Functions}, New York, Wiley.

\item
Huber, P. (2009). {\em Robust Statistics}, 2nd ed. John Wiley \&
Sons Inc.

\bibitem[Li et al., 2007]{Li}
Li, X., Hu, Z. and Wu, W. (2007). A note on the convergence of the
mean shift. {\em Pattern Recognition}, 40(6), 1756--1762.

\item
Maronna, R.A. (1976). Robust M-estimators of multivariate location
and scatter. {\em Annals of Statistics}, 4(1):51--67.

\item
Notsu, A., Komori, O. and Eguchi, S. (2014). Spontaneous clustering
via minimum gamma-divergence. {\em Neural Computation}, 26:421--448.

\item
van de Geer, S.A. (2000). {\em Empirical Processes in M-estimation:
Applications of Empirical Process Theory}, Cambridge Series in
Statistical and Probabilistic Mathematics, Cambridge University
Press.

\bibitem[Windham, 1995]{Windham}
Windham, M.~P. (1995).
\newblock Robustifying model fitting.
\newblock {\em Journal of the Royal Statistical Society, Series
  B-Methodology}, 57(3):599--609.
\end{thebibliography}
\end{document}